\documentclass[a4paper,10pt]{article}

\usepackage[dvips]{graphics,graphicx,epsfig}

\usepackage{fancybox}
\usepackage{latexsym}
\usepackage{amssymb}
\usepackage[all]{xy}
\usepackage{color}

\def\rien{\rule{0pt}{0pt}}
\begin{document}

%%%%%%%%%%%%%%%%%%%%%%%%%%%%%%%%%%
% Title

\title{Conflict Free Rule for Combining Evidences}

\author{\begin{tabular}{c}
Fr\'ed\'eric Dambreville\\
D\'el\'egation G\'en\'erale pour l'Armement, DGA/CEP/GIP/SRO\\
\noindent
16 Bis, Avenue Prieur de la C\^ote d'Or\\
Arcueil, F 94114, France\\
Web: {\tt http://www.FredericDambreville.com}\\
Email: {\tt http://email.FredericDambreville.com}
\end{tabular}
}
\maketitle
{{\bf{Abstract:}}
Recent works have investigated the problem of the conflict redistribution in the fusion rules of evidence theories.
As a consequence of these works, many new rules have been proposed.
Now, there is not a clear theoretical criterion for a choice of a rule instead another.
The present chapter proposes a new theoretically grounded rule, based on a new concept of sensor independence.
This new rule avoids the conflict redistribution, by an adaptive combination of the beliefs.
Both the logical grounds and the algorithmic implementation are considered.
}

%\begin{keywords}
%Evidence Theories, DSmT, Entropy, Optimization, Modal Logic
%\end{keywords}

\section{Introduction}
\label{DSmTb2:Dmb:Sec:1}
Recent works have underlined the limitations of the historical rule of Dempster and Shafer for fusing the information\cite{dempster2,shafer}.
The difficulty comes essentially from the conflict generation which is inherent to the rule definition.
By the way, a sequential use of the rules would result in an accumulation of the conflict, if there were not a process for removing it.
Many solutions have been proposed for managing this conflict.
The following methods are noteworthy:
\begin{itemize}
\item Constraining, relaxing  or adapting the models in order to avoid the conflict,
\item Weakening the conflicting information with the time,
\item Redistributing the conflict within the rules.
\end{itemize}
Model adaptations are of different natures.
Close to Dempster-Shafer theory, Appriou\cite{appriou} suggests to reduce the possible contradictions by a convenient setting of the problem hypotheses.
Smets\cite{smets} removes the nullity constraint on the  belief of the empty proposition (TBM);
this way, the conflict is no more a problem.
Dezert and Smarandache\cite{dezert,DSmTBook1} defined evidences on models with weakened negations (free DSmT and similar models).
By weakening or suppressing the negation, the conflict actually disappears.
The TBM of Smets and the DSmT of Dezert and Smarandache are both theoretically grounded.
TBM is axiomatically derived, while free DSmT is constructed logically\cite{DSmTb1:dmb}.
Moreover, although the DSmT keeps the nullity constraint for the empty proposition, it is possible to interpret the TBM by means of a constrained DSm model.
\\[5pt]
Avoiding the conflict by adapted models is not satisfactory however.
Indeed, there are many cases where such models appear quite artificial and not well suited to represent the real world.
Weakening the information is not satisfactory either; in many cases, the choice of a weakening criterion is rather subjective.
Experimentations\cite{martin2} have shown better results by means of rules with conflict redistributions adapted to the problem.\footnote{In fact, Dempster-Shafer rule is also a rule with redistribution of the conflict.
But in this case, the redistribution is uniform.}
Florea, Jousselme and al\cite{FloJouss} proposed recently a new family of rules which are adaptive with the conflict level.
In this case, there is an important idea: the redistribution policy is now changing automatically as a function of the conflict.
\\[5pt]
Many new rules have been proposed.
However, there is not a clear theoretical criterion for a choice of a rule instead another.
Now, these new rules, and particularly the adaptive rule of Florea and Jousselme, have uncovered a new fact: \emph{there is not a simple and permanent definition of the fusion rule for any fusion problem.}
More precisely, the structure of the fusion rule may depend on the structure of the problem.
In this chapter, we are proposing a methodology for computing fusion rules, being given a problem setting.
This methodology is logically grounded and based on a new concept of sensor independence.
As a result, the rules are obtained from a constrained convex optimization.
These computed rules cannot be derived mathematically in general.
\\[5pt]
The next section introduces evidence theories and its various viewpoint.
As a general framework for these theories, the notions of hyperpower sets and of pre-Boolean algebras are briefly reminded.
Section~\ref{DSmTb2:Dmb:Sec:3} settles a new methodology for deriving the fusion rule.
This methodology is based on an entropic notion of sensor independence.
Then, section~\ref{DSmTb2:Dmb:Sec:4} discusses about the implementations and the properties of the new rules.
Typical examples are considered.
Section~\ref{DSmTb2:Dmb:Sec:5} is more theoretical and exposes the logical fundaments of our methodology.
At last, section~\ref{DSmTb2:Dmb:Sec:6} concludes.
\section{Viewpoints in evidence theories}
\label{DSmTb2:Dmb:Sec:2}
In this section, we are discussing about theories for combining evidences expressed by belief functions.
Since pre-Boolean algebra is a common framework for all these theories, in particular as a generalization of sets and hyperpower sets, we are now introducing briefly this notion.
\subsection{Pre-Boolean algebra}
The theory of Dezert and Smarandache is based on the fundamental notion of pre-Boolean algebra, or hyperpower sets.
These algebra will describe the logical modeling of the knowledge.
This chapter is not dedicated to a thorough exposition of the theory of pre-Boolean algebra.
The reader should refer to the chapter~\ref{Dmb:PremierChapitre} of this book for a precise theoretical definition.
Now, the present section will introduce these notions qualitatively, and some typical examples will be provided.
\subsubsection{General principle}
Subsequently, the conjunction and disjunction are denoted $\wedge$ and $\vee$.
The negation, when used, is denoted $\neg$.
The empty set is denoted $\bot$ while the tautology, or full ignorance, is denoted $\top$.
Notice that these notations are not the most classical in the domain of evidence theories.
Typically, $\cap,\cup,\Theta\setminus\cdot,\emptyset,\Theta$ are used instead of $\wedge,\vee,\neg,\bot,\top$.
However, $\wedge,\vee,\neg,\bot,\top$ are notations widely used in logics and Boolean algebra.
Since the connexions are important between these theories, we will use the logical notations in general.
\paragraph{Definition.} 
A pre-Boolean algebra could be seen as a subset of a Boolean algebra which is stable for the conjunction and the disjunction.
As a consequence, a pre-Boolean algebra together with the two operators, conjunction and disjunction, is an algebraic structure.
\\[5pt]
This algebraic structure has the same properties than a Boolean algebra, except that it does not implement explicitly the notion of negation.
In particular, the following properties are provided by the pre-Boolean algebra for the binary operators:
\begin{description}
\item[\rien$\quad$\emph{Commutativity.}]$\phi\wedge\psi=\psi\wedge\phi$ and $\phi\vee\psi=\psi\vee\phi$\,,
\item[\rien$\quad$\emph{Associativity.}]$\phi\wedge(\psi\wedge\eta)=(\phi\wedge\psi)\wedge\eta$ and $\phi\vee(\psi\vee\eta)=(\phi\vee\psi)\vee\eta$\,,
\item[\rien$\quad$\emph{Distributivity.}]$\phi\wedge(\psi\vee\eta)=(\phi\wedge\psi)\vee(\phi\wedge\eta)$ and $\phi\vee(\psi\wedge\eta)=(\phi\vee\psi)\wedge(\phi\vee\eta)$\,,
\item[\rien$\quad$\emph{Idempotence.}]$\phi\wedge\phi=\phi$ and $\phi\vee\phi=\phi$\,,
\item[\rien$\quad$\emph{Neutral sup/sub-elements.}]$\phi\wedge(\phi\vee\psi)=\phi$ and $\phi\vee(\phi\wedge\psi)=\phi$\,,
\end{description}
for any $\phi,\psi,\eta$ in the pre-Boolean algebra.
\subsubsection{Example}
\begin{figure}[ht!]
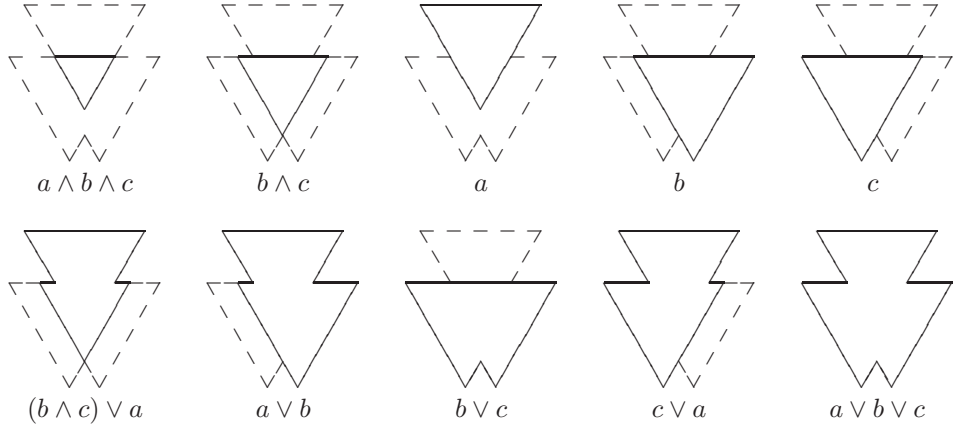

\centering
\scalebox{1}{\begin{tabular}{@{}ccccc@{}}
\xy
(0,0)*+{\rien},
%\green
(2,0)="a1",
(18,0)="b1",
(10,-13.9)="c1",
(20,-6.9)="a2",
(12,-20.8)="b2",
(6,-6.9)="c2",
(8,-20.8)="a3",
(0,-6.9)="b3",
(14,-6.9)="c3",
(10,-17.35)="c1b",
(4,-6.9)="c2b",
(16,-6.9)="c3b",
\ar @{--}"a1";"b1",
\ar @{--}"b1";"c3",
\ar @{-}"c3";"c2",
\ar @{--}"c2";"a1",
\ar @{--}"c3";"a2",
\ar @{--}"a2";"b2",
\ar @{--}"b2";"c1b",
\ar @{--}"c1b";"a3",
\ar @{--}"a3";"b3",
\ar @{--}"b3";"c2",

\ar @{-}"c1";"c3",
\ar @{-}"c2";"c1",
\POS "c1",
(20,-20.8)*+{\rien}
\endxy
&
\xy
(0,0)*+{\rien},
%\green
(2,0)="a1",
(18,0)="b1",
(10,-13.9)="c1",
(20,-6.9)="a2",
(12,-20.8)="b2",
(6,-6.9)="c2",
(8,-20.8)="a3",
(0,-6.9)="b3",
(14,-6.9)="c3",
(10,-17.35)="c1b",
(4,-6.9)="c2b",
(16,-6.9)="c3b",
\ar @{--}"a1";"b1",
\ar @{--}"b1";"c3",
\ar @{--}"c2";"a1",
\ar @{--}"c3b";"a2",
\ar @{--}"a2";"b2",
\ar @{--}"b2";"c1b",
\ar @{--}"c1b";"a3",
\ar @{--}"a3";"b3",
\ar @{--}"b3";"c2b",
\ar @{-}"c2b";"c3b",
\ar @{-}"c3b";"c1b",
\ar @{-}"c1b";"c2b",
\POS "c2b",
(20,-20.8)*+{\rien}
\endxy
&
\xy
(0,0)*+{\rien},
%\green
(2,0)="a1",
(18,0)="b1",
(10,-13.9)="c1",
(20,-6.9)="a2",
(12,-20.8)="b2",
(6,-6.9)="c2",
(8,-20.8)="a3",
(0,-6.9)="b3",
(14,-6.9)="c3",
(10,-17.35)="c1b",
(4,-6.9)="c2b",
(16,-6.9)="c3b",
\ar @{-}"a1";"b1",
\ar @{-}"b1";"c1",
\ar @{-}"c1";"a1",
\ar @{--}"c2";"b3",
\ar @{--}"b3";"a3",
\ar @{--}"a3";"c1b",
\ar @{--}"c1b";"b2",
\ar @{--}"b2";"a2",
\ar @{--}"a2";"c3",
\POS "c3",
(20,-20.8)*+{\rien}
\endxy
&
\xy
(0,0)*+{\rien},
%\green
(2,0)="a1",
(18,0)="b1",
(10,-13.9)="c1",
(20,-6.9)="a2",
(12,-20.8)="b2",
(6,-6.9)="c2",
(8,-20.8)="a3",
(0,-6.9)="b3",
(14,-6.9)="c3",
(10,-17.35)="c1b",
(4,-6.9)="c2b",
(16,-6.9)="c3b",
\ar @{-}"b2";"c2b",
\ar @{-}"c2b";"a2",
\ar @{-}"a2";"b2",

\ar @{--}"c2";"a1",
\ar @{--}"a1";"b1",
\ar @{--}"b1";"c3",
\ar @{--}"c2b";"b3",
\ar @{--}"b3";"a3",
\ar @{--}"a3";"c1b",
\POS "c1b",
(20,-20.8)*+{\rien}
\endxy
&
\xy
(0,0)*+{\rien},
%\green
(2,0)="a1",
(18,0)="b1",
(10,-13.9)="c1",
(20,-6.9)="a2",
(12,-20.8)="b2",
(6,-6.9)="c2",
(8,-20.8)="a3",
(0,-6.9)="b3",
(14,-6.9)="c3",
(10,-17.35)="c1b",
(4,-6.9)="c2b",
(16,-6.9)="c3b",
\ar @{-}"c3b";"a3",
\ar @{-}"a3";"b3",
\ar @{-}"b3";"c3b",
\ar @{--}"c2";"a1",
\ar @{--}"a1";"b1",
\ar @{--}"b1";"c3",
\ar @{--}"c3b";"a2",
\ar @{--}"a2";"b2",
\ar @{--}"b2";"c1b",
\POS "c1b",
(20,-20.8)*+{\rien}
\endxy
\\
$a\wedge b\wedge c$&$b\wedge c$&$a$&$b$&$c$
\vspace{3pt}
\\
\xy
(0,0)*+{\rien},
%\green
(2,0)="a1",
(18,0)="b1",
(10,-13.9)="c1",
(20,-6.9)="a2",
(12,-20.8)="b2",
(6,-6.9)="c2",
(8,-20.8)="a3",
(0,-6.9)="b3",
(14,-6.9)="c3",
(10,-17.35)="c1b",
(4,-6.9)="c2b",
(16,-6.9)="c3b",
\ar @{-}"c2";"a1",
\ar @{-}"a1";"b1",
\ar @{-}"b1";"c3",
\ar @{-}"c3";"c3b",
\ar @{-}"c3b";"c1b",
\ar @{-}"c1b";"c2b",
\ar @{-}"c2b";"c2",
\ar @{--}"c1b";"a3",
\ar @{--}"a3";"b3",
\ar @{--}"b3";"c2b",
\ar @{--}"c1b";"b2",
\ar @{--}"b2";"a2",
\ar @{--}"a2";"c3b",
\POS "c3b",
(20,-20.8)*+{\rien}
\endxy
&
\xy
(0,0)*+{\rien},
%\green
(2,0)="a1",
(18,0)="b1",
(10,-13.9)="c1",
(20,-6.9)="a2",
(12,-20.8)="b2",
(6,-6.9)="c2",
(8,-20.8)="a3",
(0,-6.9)="b3",
(14,-6.9)="c3",
(10,-17.35)="c1b",
(4,-6.9)="c2b",
(16,-6.9)="c3b",
\ar @{-}"c2";"a1",
\ar @{-}"a1";"b1",
\ar @{-}"b1";"c3",
\ar @{-}"c3";"a2",
\ar @{-}"a2";"b2",
\ar @{-}"b2";"c2b",
\ar @{-}"c2b";"c2",
\ar @{--}"c1b";"a3",
\ar @{--}"a3";"b3",
\ar @{--}"b3";"c2b",
\POS "c3",
(20,-20.8)*+{\rien}
\endxy
&
\xy
(0,0)*+{\rien},
%\green
(2,0)="a1",
(18,0)="b1",
(10,-13.9)="c1",
(20,-6.9)="a2",
(12,-20.8)="b2",
(6,-6.9)="c2",
(8,-20.8)="a3",
(0,-6.9)="b3",
(14,-6.9)="c3",
(10,-17.35)="c1b",
(4,-6.9)="c2b",
(16,-6.9)="c3b",
\ar @{-}"a2";"b2",
\ar @{-}"b2";"c1b",
\ar @{-}"c1b";"a3",
\ar @{-}"a3";"b3",
\ar @{-}"b3";"a2",
\ar @{--}"c2";"a1",
\ar @{--}"a1";"b1",
\ar @{--}"b1";"c3",
\POS "c3",
(20,-20.8)*+{\rien}
\endxy
&
\xy
(0,0)*+{\rien},
%\green
(2,0)="a1",
(18,0)="b1",
(10,-13.9)="c1",
(20,-6.9)="a2",
(12,-20.8)="b2",
(6,-6.9)="c2",
(8,-20.8)="a3",
(0,-6.9)="b3",
(14,-6.9)="c3",
(10,-17.35)="c1b",
(4,-6.9)="c2b",
(16,-6.9)="c3b",
\ar @{-}"a1";"b1",
\ar @{-}"b1";"c3",
\ar @{-}"c3";"c3b",
\ar @{-}"c3b";"a3",
\ar @{-}"a3";"b3",
\ar @{-}"b3";"c2",
\ar @{-}"c2";"a1",
\ar @{--}"c3b";"a2",
\ar @{--}"a2";"b2",
\ar @{--}"b2";"c1b",
\POS "c1b",
(20,-20.8)*+{\rien}
\endxy
&
\xy
(0,0)*+{\rien},
%\green
(2,0)="a1",
(18,0)="b1",
(10,-13.9)="c1",
(20,-6.9)="a2",
(12,-20.8)="b2",
(6,-6.9)="c2",
(8,-20.8)="a3",
(0,-6.9)="b3",
(14,-6.9)="c3",
(10,-17.35)="c1b",
(4,-6.9)="c2b",
(16,-6.9)="c3b",
\ar @{-}"a1";"b1",
\ar @{-}"b1";"c3",
\ar @{-}"c3";"a2",
\ar @{-}"a2";"b2",
\ar @{-}"b2";"c1b",
\ar @{-}"c1b";"a3",
\ar @{-}"a3";"b3",
\ar @{-}"b3";"c2",
\ar @{-}"c2";"a1",
\POS "c3",
(20,-20.8)*+{\rien}
\endxy
\\
$(b\wedge c)\vee a$&$a\vee b$&$b\vee c$&$c\vee a$&$a\vee b\vee c$
\end{tabular}}
\caption{Pre-Boolean algebra $<a,b,c>_\Gamma$; ($\bot$ and $\top$ are omitted)}
\label{dsmtb2:dmb:v2:fig3}
\end{figure}

\paragraph{Free pre-Boolean algebra.}
Let $a,b,c$ be three atomic propositions.
Consider the free Boolean algebra $\mathcal{B}(a,b,c)$ generated by $a,b,c$\,:
$$\mbox{$
\mathcal{B}(a,b,c)=\left\{\left.\bigvee_{(\alpha,\beta,\gamma)\in A}\bigl(\alpha\wedge\beta\wedge\gamma\bigr)\right/ A\subset\{a,\neg a\}\times\{b,\neg b\}\times\{c,\neg c\}\right\}\;.
$}$$
It is well known that $\mathcal{B}(a,b,c)$ contains $2^{2^3}=256$ elements.
\\[3pt]
The free pre-Boolean algebra generated by the propositions $a,b,c$ is the smaller subset of $\mathcal{B}(a,b,c)$ containing $a,b,c$ and stable for $\wedge$ and $\vee$.
This set, denoted $<a,b,c>$, is defined extensionally by:
$$\begin{array}{@{}l@{}}
<a,b,c>=\{\bot,
a\wedge b\wedge c,
a\wedge b, a\wedge c, b\wedge c,
a\wedge (b\vee c), b\wedge (a\vee c), c\wedge (a\vee b),
a,b,c,
\vspace{3pt}\\
\rien\hspace{15pt}
(a\wedge b)\vee(a\wedge c)\vee(b\wedge c),
(a\wedge b)\vee c,(a\wedge c)\vee b,(b\wedge c)\vee a,
a\vee b, a\vee c, b\vee c,
a\vee b\vee c,
\top\}
\end{array}$$
It is easily verified that a conjunctive or disjunctive combination of propositions of $<a,b,c>$ is still a proposition of $<a,b,c>$.
For example:
$$
\bigl(a\wedge (b\vee c)\bigr)\vee b=\bigl((a\wedge b)\vee (a\wedge c)\bigr)\vee b=(a\wedge c) \vee b \;.
$$
Moreover, $<a,b,c>$ is obviously the smallest set, which is stable for $\wedge$ and $\vee$.
In particular, it is noticed that $\bot,\top\in<a,b,c>$ since $\bot=\bigvee_{\alpha\in\emptyset}\alpha$ and $\top=\bigwedge_{\alpha\in\emptyset}\alpha$\,.
\\[5pt]
The free pre-Boolean algebra $<a,b,c>$ is also called \emph{hyperpower set} generated by $a,b,c$\,.
It is also denoted $D^\Theta$, where $\Theta=\{a,b,c\}$\,.
Notice that the tautology $\top$ is often excluded from the definition of the hyperpower set\cite{DSmTBook1}\,.
By the way, Dambreville excluded both $\bot$ and $\top$ from a previous definition~\cite{DSmTb1:dmb}.
These differences have a quite limited impact, when considering the free DSmT.
Whatever, it is generally assumed that $a\vee b\vee c=\top$\,;
\emph{but this is an additional hypothesis}.
\paragraph{A Boolean algebra is a constrained pre-Boolean algebra.}
A Boolean algebra is a subset of itself and is stable for $\wedge$ and $\vee$.
Thus, it is a pre-Boolean algebra.
Now, we will see on an example that a set could be seen as an hyperpower set which has been constrained by logical constraints.
Since a Boolean algebra could be considered as a set, this result implies more generally that a Boolean algebra could be obtained by constraining a free pre-Boolean algebra.
\\[5pt]
Denote $\Theta=\{a,b,c\}$\,.
Consider the Boolean algebra $\mathcal{P}(\Theta)$ related to the set operators $\cap,\cup,\Theta\setminus\cdot$ and neutral elements $\emptyset,\Theta$.
This Boolean algebra is extensionally defined by:
$$
\mathcal{P}(\Theta)=\bigl\{
\emptyset,
a,b,c,
\{a,b\}, \{a,c\}, \{b,c\},
\Theta
\bigr\}\;.
$$
Now, consider the hyperpower set $<a,b,c>$ and apply to it the constraints:
$$
\Gamma=\big\{
a\wedge b=b\wedge c=a\wedge c=\bot
\;,\ a\vee b\vee c=\top
\bigr\}\;.
$$
It is then derived:
$$\left\{\begin{array}{l@{}}
a\wedge b\wedge c=a\wedge (b\vee c)=b\wedge (a\vee c)=c\wedge (a\vee b)=(a\wedge b)\vee(a\wedge c)\vee(b\wedge c)=\bot\;,
\\
(a\wedge b)\vee c=c
\;, 
(a\wedge c)\vee b=b
\;, 
(b\wedge c)\vee a=a
\;,
\\
a\vee b\vee c=\top\;.
\end{array}\right.$$
Denoting $<a,b,c>_\Gamma$ the resulting constrained pre-Boolean algebra, it comes:
$$
<a,b,c>_\Gamma=\{\bot,a,b,c,a\vee b, a\vee c, b\vee c, \top\}\;.
$$
Then, $<a,b,c>_\Gamma$ contains exactly the same number of elements than $\mathcal{P}(\Theta)$.
More precisely, by the Boolean properties of $\wedge$ and $\vee$, it is clear that $<a,b,c>_\Gamma$ and $\mathcal{P}(\Theta)$ are isomorph as pre-Boolean algebra.
While $<a,b,c>_\Gamma$ does not define the negation explicitly, this isomorphism shows that the negation is implicitly defined in $<a,b,c>_\Gamma$\,.
In fact, the negation of $<a,b,c>_\Gamma$ has been built by the constraints.
This is an important property of pre-Boolean algebra:
\begin{quote}\emph{
The constraints put on a free pre-Boolean algebra partially characterize the negation operator.
}\end{quote}
As a consequence, there is a partial definition of the negation in a pre-Boolean algebra.
This negation is entirely undefined in an hyperpower set and is entirely defined in a set.
But there are many intermediate cases.
\paragraph{Example of constrained pre-Boolean algebra.}
Let $\Theta=\{a,b,c\}$ be a set of atomic propositions and $\Gamma=\{a\wedge b=a\wedge c\}$ be a set of constraints.
By propagating the constraints, it is obtained:
$$
a\wedge b=a\wedge c=(a\wedge b)\vee(a\wedge c)=(a\wedge b)\wedge(a\wedge c)=a\wedge b\wedge c\,.
$$
Consequently:
$$\begin{array}{@{}l@{}}
a\wedge (b\vee c)=a\wedge b\wedge c\,,\ 
(a\wedge b)\vee c=c\,,\ 
(a\wedge c)\vee b=b\,.
\vspace{3pt}\\
b\wedge (a\vee c)=c\wedge (a\vee b)=(a\wedge b)\vee(a\wedge c)\vee(b\wedge c)=b\wedge c\,.
\end{array}$$
At last, the constrained pre-Boolean algebra is extensionally defined by:
$$
<a,b,c>_\Gamma=\{\bot,
a\wedge b\wedge c, b\wedge c,
a,b,c,
(b\wedge c)\vee a,
a\vee b, a\vee c, b\vee c,
a\vee b\vee c,
\top\}
$$
This configuration is modeled in figure~\ref{dsmtb2:dmb:v2:fig3}.
This model ensures that the propagation of the constraints is complete in the definition of $<a,b,c>_\Gamma$.
\subsubsection{Notations}
Let $\Theta$ be a set of atomic propositions.
The free pre-Boolean algebra generated by $\Theta$ is denoted $<\Theta>$\,.
\\[5pt]
Now, let $\Gamma$ be a set of constraints over the propositions of $<\Theta>$.
The pre-Boolean algebra generated by $\Theta$ and constrained by $\Gamma$ is denoted $<\Theta>_\Gamma$\,.
Of course, it comes $<\Theta>_\emptyset=<\Theta>$ (the pre-Boolean algebra generated by $\Theta$ and constrained by an empty $\Gamma$ is an hyperpower set).
\\[5pt]
A proposition $\phi$ is a subproposition of proposition $\psi$ if and only if $\phi\wedge\psi=\phi$\,;
subsequently, the property $\phi\wedge\psi=\phi$ is also denoted $\phi\subset\psi$\,.
\subsection{Belief}
It is now given a pre-Boolean algebra $<\Theta>_\Gamma$ as a logical framework for the knowledge representation.
The theories of evidence also implement a belief on each logical proposition.
This belief contains both an imprecision and an uncertainty information.
The following sections consider two main styles for implementing the belief.
In the DSmT and DST (Dempster Shafer Theory)~\cite{shafer}, the belief over the empty proposition is always zero.
In the TBM (Transferable Belief Model)~\cite{smets}, the belief over the empty proposition may be non zero.
These viewpoints are related to two slightly different logical interpretations, as stated in section~\ref{DSmTb2:Dmb:Sec:5}.
\subsubsection{DSmT and DST}
DSmT defines the notion of belief in a same way than DST.
The only difference is that DST works on a set, while DSmT works on any pre-Boolean algebra.
Fundamental differences will also arise, when defining the fusion of the information (section~\ref{DSmTb2:Dmb:CH2:section:FusionRule}).
\paragraph{Basic Belief Assignment.}
A basic belief assignment (bba) to the pre-Boolean algebra $<\Theta>_\Gamma$ is a real valued function defined over $<\Theta>_\Gamma$ such that:
$$
\sum_{\phi\in<\Theta>_\Gamma}m(\phi)=1
\ ,\quad m(\bot)=0\quad\mbox{and}\quad m\ge0
\;.
$$
Typically, $m$ represents the knowledge of an expert or of a sensor.
By hypothesizing \mbox{$m(\bot)=0$}, the DSmT assumes the \emph{coherence} of the information.
\\[5pt]
The bba is a belief density, describing the information intrinsic to the propositions.
The \emph{full} belief of a proposition is thus the compilation of the bba of its sub-propositions.
\paragraph{Belief function.}
The belief function $\mathrm{Bel}$ related to a bba $m$ is defined by:
\begin{equation}
\label{f2k5:DSmTcont:Eq:1}
\forall\phi\in<\Theta>_\Gamma,\,\mathrm{Bel}(\phi)=\sum_{\psi\in<\Theta>_\Gamma:\psi\subset\phi}m(\psi)\;.
\end{equation}
It is generally considered that $\mathrm{Bel}\left(\bigvee_{\theta\in\Theta}\theta\right)=1$\,, which means that $\Theta$ matches all possible information.
\subsubsection{TBM and TBM-like bba}
Like the DST, the TBM works on a set.
However, in the TBM interpretation the belief put on the empty set is not necessarily zeroed.
It is also possible to mix this hypothesis with a pre-Boolean modeling, as follows.
\paragraph{TBM-like Basic Belief Assignment.}
A basic belief assignment to the pre-Boolean algebra $<\Theta>_\Gamma$ is a real valued function defined over $<\Theta>_\Gamma$ such that:
$$
\sum_{\phi\in<\Theta>_\Gamma}m(\phi)=1
\quad\mbox{and}\quad m\ge0
\;.
$$
By removing the hypothesis \mbox{$m(\bot)=0$}, the coherence of the information is canceled.
The coherence and non-coherence hypotheses have a logical interpretation, as explained in section~\ref{DSmTb2:Dmb:Sec:5}.
\\[5pt]
In fact, it is also possible to simulate the TBM (and TBM-like models) by means of the DSmT (with the coherence hypothesis).
The idea is to simulate the empty set of TBM by the pre-Boolean proposition $\bigwedge_{\theta\in\Theta}\theta$\,.
This result, although simple, is outside the scope of this chapter and will not be developed further.
%\\[5pt]
To end with this subsection, it is noticed that Smets proposes a slightly different definition of the belief function by excluding the belief of the empty set.
Smets belief function will be denoted and defined by:
$$
\mathrm{Bel}_S(\phi)= \sum_{\psi\in<\Theta>_\Gamma:\bot\ne\psi\subset\phi}m(\psi)\;.
$$
This truncated belief function is not used subsequently, since we work essentially on bba and on the full belief function $\mathrm{Bel}$ as defined in~(\ref{f2k5:DSmTcont:Eq:1}).
\subsection{Fusion rules}
\label{DSmTb2:Dmb:CH2:section:FusionRule}
The main contribution of evidence theories consists in their fusion rules.
It is assumed then that two or more sources of information are providing a
viewpoint about the universe.
These viewpoints are described by specific bbas for each sensor.
The question then is to make a unique representation of the information,
\emph{i.e.} a unique bba, from these several bbas.
Several rules for fusing such information have been elaborated.
\\[5pt]
There are essentially two kinds of rules.
The first kind avoids any conflict redistribution.
The theorists generally agree then on a unique fusion rule, the conjunctive
rule (without redistribution).
Two models avoid the conflict redistribution: the transferable belief model of
Smets and the free DSmT.
In both theories, a strong constraint is put on the model.
TBM puts non zero weights on the empty set, while free DSmT removes the negation from the model.
In many cases however, these hypotheses are too restrictive.
\\[5pt]
When the conflict is taken into account and is redistributed, many possible
rules have been proposed.
No model restriction is needed anymore, but it is difficult to decide for a
definitive fusion rule.
\\[5pt]
The following sections introduce shortly these various concept of rules.
\subsubsection{Fusion rule in free DSmT and similar models.}
Free DSmT is defined on an hyperpower set.
A fundamental property of an hyperpower set is that the empty proposition cannot be generated from non empty propositions.
More generally, a pre-Boolean algebra $<\Theta>_\Gamma$, where the constraints in $\Gamma$ do not generate $\bot$, will also satisfy such property:
\begin{equation}\label{dmb:dsmtb2:ch2:insulation:eq:1}
\phi,\psi\in<\Theta>\setminus\{\bot\} \Longrightarrow \phi\wedge\psi\in<\Theta>\setminus\{\bot\}\,.
\end{equation}
This property will be called an \emph{insulation} property.
\\[5pt]
Assume now a pre-Boolean algebra $<\Theta>_\Gamma$ satisfying (\ref{dmb:dsmtb2:ch2:insulation:eq:1}).
Then, two bbas $m_1$ and $m_2$ over $<\Theta>_\Gamma$ will be
fused into a bba $m_1\oplus m_2$ as follows: 
\begin{equation}\label{dmb:dsmtb2:ch2:rule:eq:2}
\forall\phi\in<\Theta>_\Gamma\,,\;m_1\oplus m_2(\phi)=\sum_{\psi\wedge\eta=\phi}m_1(\psi) m_2(\eta)\;.
\end{equation}
This definition is compatible with the constraint $m_1\oplus m_2(\bot)=0$ of DSmT, since it comes by the insulation property:
$$
m_1(\psi)\ne0\mbox{ and }m_2(\eta)\ne0\mbox{ \ implies \ }\psi\wedge\eta\in<\Theta>\setminus\{\bot\}\,.
$$
The insulation property is often a too strong hypothesis for many problems.
The TBM viewpoint will not request such structure constraints.
But as a consequence, the coherence property of the bba will be removed.
\subsubsection{Fusion rule for TBM-like bbas}
In the TBM paradigm, two bbas $m_1$ and $m_2$ over $<\Theta>_\Gamma$ will be
fused into a bba $m_1\oplus m_2$ as follows: 
\begin{equation}\label{dmb:dsmtb2:ch2:rule:eq:2:bis}
\forall\phi\in<\Theta>_\Gamma\,,\;m_1\oplus m_2(\phi)=\sum_{\psi\wedge\eta=\phi}m_1(\psi) m_2(\eta)\;.
\end{equation}
There is no particular restriction on the choice of $<\Theta>_\Gamma$ in this case.
It is for example possible that the model contains two non empty propositions $\psi$ and $\eta$ such that $\psi\wedge\eta=\bot$\,.
Assuming that the initial bbas $m_1$ and $m_2$ are such that $m_1(\psi)>0$ and $m_2(\eta)>0$, it comes from the definition that $m_1\oplus m_2(\bot)>0$\,.
But the rule is still compatible with the TBM paradigm, since then the \emph{coherence} constraint $m_1\oplus m_2(\bot)=0$ is removed.
By the way, removing this constraint is not satisfactory in many cases.
In particular, it is well known that the weight of the contradiction may increase up to $1$ by iterating the fusion stages.
\subsubsection{General case}
While the fusion rule is clearly defined by~(\ref{dmb:dsmtb2:ch2:rule:eq:2}) for models avoiding the conflict, there are many possible rules when this conflict has to be redistributed.
Typically, the rule could be defined in two steps.
First, compute the conjunctive function $\mu$ of $m_1$ and $m_2$ by:
$$
\forall\phi\in<\Theta>_\Gamma\,,\;\mu(\phi)=\sum_{\psi\wedge\eta=\phi}m_1(\psi) m_2(\eta)\;.
$$
The function $\mu$ is like the fusion rule in the TBM paradigm.
It cannot be used directly, since $\mu(\bot)$ have to be redistributed when $\mu(\bot)>0$\,.
Redistributing the conflict means:
\begin{itemize}
\item Constructing a function $\rho$ on $<\Theta>_\Gamma$ such that:
$$
\rho(\bot)=0\;,\quad\sum_{\phi\in<\Theta>_\Gamma}\rho(\phi)=1\quad\mbox{and}\quad\rho\ge0\;,
$$
\item Derive the bba $m_1\oplus m_2$ by:
$$
m_1\oplus m_2(\phi)=\mu(\phi)+\rho(\phi)\mu(\bot)\;.
$$
\end{itemize}
There are many possible rules deduced from the redistribution principle.
Moreover, the redistribution may be dependent to a local conflict, like the PCR rules\cite{DezSm:3,martin2}.
It is also noticed that some authors\cite{FloJouss} allows negative redistributions by removing the constraint $\rho\ge0$\,.
% A VERIFIER!
These new rules are as well legitimate and interesting, but by allowing negative redistributions, the criterion for defining rules is again weakened.
The question now is \emph{how to decided for a rule or another?}
This choice is certainly dependent of the structure of the fusion problem.
Actually, Florea, Jousselme and al\cite{FloJouss} proposed a rule adaptive with the conflict level.
More generally, it is foreseeable that a fusion rule should be defined or computed specifically for a given fusion problem.
\\[5pt]
In the next sections, we will derive logically a new ruling method, which avoids the conflict redistribution by exploiting a new concept of independence of the sensors.
The new rules will be essentially computed from an entropic optimization problem.
This problem may be unsolvable, which will imply a rejection of the hypotheses (too high conflict between the sources).
Otherwise, it will adapt the belief dispatching in a much more flexible way than the usual conjunctive function $\mu$.
\section{Entropic approach of the rule definition}
\label{DSmTb2:Dmb:Sec:3}
To begin with this new rule concept, we will directly settle the concrete optimization principles of our method.
The logical justifications will come later, in section~\ref{DSmTb2:Dmb:Sec:5}.
\subsection{Independent sources and entropy}
Actually, the idea is not completely new, and Dezert used it in order to give a first justification to the free DSmT~\cite{dezert}.
More precisely, the free rule could be rewritten:
$$
\forall\phi\in<\Theta>_\Gamma\,,\;m_1\oplus m_2(\phi)=\sum_{\psi\wedge\eta=\phi}f_o(\psi,\eta)\;,
$$
where:
\begin{equation}\label{dmb:dsmtb2:ch2:indepf:eq:3}
f_o(\psi,\eta)=m_1(\psi) m_2(\eta)\;.
\end{equation}
If we are interpreting $m_i$ as a kind of probability, the relation~(\ref{dmb:dsmtb2:ch2:indepf:eq:3}) is like the probabilistic independence, where $f_o(\psi,\eta)$ is a kind of joint probability.
Section~\ref{DSmTb2:Dmb:Sec:5} will clarify this probabilistic viewpoint.
Now, there is a link between the notion of probabilistic independence and the notion of entropy, which is often forgotten.
The law $f_o(\psi,\eta)=m_1(\psi) m_2(\eta)$ is a maximizer of the entropy, with respect to the constraint of marginalization:
\begin{equation}\label{dmb:dsmtb2:ch2:indepf:eq:4}
%\left\{
\begin{array}{@{\;}l@{}}\displaystyle
f_o\in\arg\max_f-\sum_{\psi,\eta}f(\psi,\eta)\ln f(\psi,\eta)
\vspace{5pt}\\\mbox{under constraints: }
\vspace{4pt}\\\displaystyle
\rien\hspace{40pt}
f\ge0\;,\quad
\sum_{\psi}f(\psi,\eta)=m_2(\eta)\quad\mbox{and}\quad
\sum_{\eta}f(\psi,\eta)=m_1(\psi)\;.
\end{array}
%\right.
\end{equation}
This is actually how Dezert derived the conjunctive rule of free DSmT~\cite{dezert}, although he did not make an explicit mention to the probability theory.
Now, the equation~(\ref{dmb:dsmtb2:ch2:indepf:eq:4}) has a particular interpretation in the paradigm of information theory:
$f_o$ is the law which contains the maximum of information, owing to the fact that its marginals are $m_1$ and $m_2$.
By the way, independent sources of information should provide the maximum of information, so that the maximization of entropy appears as the good way to characterize independent sources.
When the constraints are just the marginalizations, the solution to this maximization is the independence relation $f_o(\psi,\eta)=m_1(\psi) m_2(\eta)$.
In Bayesian theory particularly, the marginalizations are generally the only constraints, and the notion of independent sources of information reduces to the notion of independent propositions.
But in the case of evidence theories, \emph{there is the problem of the conflict}, which adds constraints.
\subsection{Definition of a new rule for the DSmT}
Let be defined a pre-Boolean algebra $<\Theta>_\Gamma$, constituting the logical framework of the information.
Let be defined two bbas $m_1$ and $m_2$ over $<\Theta>_\Gamma$.
The bbas are assumed to be \emph{coherent}, so that $m_1(\bot)=m_2(\bot)=0$\,.
Then the fusion of $m_1$ and $m_2$ is the bba $m_1\oplus m_2$ defined by:
\begin{equation}\label{dmb:dsmtb2:ch2:rule:eq:5}
\begin{array}{@{}l@{}}\displaystyle
\forall\phi\in<\Theta>_\Gamma\,,\;m_1\oplus m_2(\phi)=\sum_{\psi\wedge\eta=\phi}f_o(\psi,\eta)\;,
\vspace{5pt}\\\mbox{where:}
\vspace{5pt}\\
\rien\hspace{30pt}
%\left\{
\begin{array}{@{\;}l@{}}\displaystyle
f_o\in\arg\max_f-\sum_{\psi,\eta}f(\psi,\eta)\ln f(\psi,\eta)
\vspace{5pt}\\\mbox{under constraints: }
\vspace{4pt}\\\displaystyle
\rien\hspace{40pt}
f\ge0\;,\quad
\sum_{\psi}f(\psi,\eta)=m_2(\eta)\;,\quad
\sum_{\eta}f(\psi,\eta)=m_1(\psi)\;,
\vspace{2pt}\\\displaystyle
\rien\hspace{40pt}
\mbox{and}\quad\forall\psi,\eta\in<\Theta>_\Gamma\,,\;
\psi\wedge\eta=\bot
\Longrightarrow f(\psi,\eta)=0\,.
\end{array}
%\right.
\end{array}
\end{equation}
This rule will be named \emph{Entropy Maximizing Rule} (EMR).
\\[5pt]
\emph{Corollary of the definition.}
The fused bba is compatible with the coherence constraint \mbox{$m_1\oplus m_2(\bot)=0$} of DSmT.
\\[5pt]
Proof is immediate owing to the constraints $\psi\wedge\eta=\bot
\Longrightarrow f(\psi,\eta)=0$ in the optimization.
\subsection{Feasibility of the rule}
The rule is feasible when there is a solution to the optimization.
The feasibility is obtained as soon there is a solution to the constraints.
\paragraph{Definition.}
The fused bba $m_1\oplus m_2$ is defined if and only if there exists a function $f$ such that:
\begin{equation}\label{dmb:dsmtb2:ch2:rule:eq:6}
\begin{array}{@{\;}l@{}}\displaystyle
f\ge0\;,\quad
\sum_{\psi}f(\psi,\eta)=m_2(\eta)\;,\quad
\sum_{\eta}f(\psi,\eta)=m_1(\psi)\;,
\vspace{2pt}\\\displaystyle
\rien\hspace{40pt}
\mbox{and}\quad\forall\psi,\eta\in<\Theta>_\Gamma\,,\;
\psi\wedge\eta=\bot
\Longrightarrow f(\psi,\eta)=0\,.
\end{array}
\end{equation}
In next section, it will be shown on examples that the fusion is not always feasible.
Actually, the infeasibility of the rule is a consequence of fundamental incompatibilities of the information.
\subsection{Generalizations}
\subsubsection{Fusion of $N$ bbas}
It will be seen that the fusion rule defined previously is not associative.
This means that the sources of information do not have the same \emph{weight} in a sequential fusion.
However, when it is needed to fuse $N$ sources of information simultaneously, the fusion method has then to be generalized to $N$ bbas.
The subsequent definition makes this generalization.
\paragraph{$N$-ary rule.}
Let be defined a pre-Boolean algebra $<\Theta>_\Gamma$, constituting the logical framework of the information.
Let be defined $N$ coherent bbas $m_i|_{1\le i\le N}$ over $<\Theta>_\Gamma$.
Then the fusion of $m_i|_{1\le i\le N}$ is the bba $\oplus[m_i|1\le i\le N]$ defined by:
\begin{equation}\label{dmb:dsmtb2:ch2:rule:eq:7}
\begin{array}{@{}l@{}}\displaystyle
\forall\phi\in<\Theta>_\Gamma\,,\;\oplus[m_i|1\le i\le N](\phi)=\sum_{\bigwedge_{i=1}^N\psi_i=\phi}
f_o(\psi_i|1\le i\le N)\;,
\vspace{5pt}\\\mbox{where:}
\vspace{5pt}\\
\rien\hspace{30pt}
%\left\{
\begin{array}{@{\;}l@{}}\displaystyle
f_o\in\arg\max_f-\sum_{\psi}f(\psi_i|1\le i\le N)\ln f(\psi_i|1\le i\le N)
\vspace{5pt}\\\mbox{under constraints: }
\vspace{4pt}\\\displaystyle
\rien\hspace{40pt}
f\ge0\;,\quad
\forall i,\;\sum_{\psi_j|j\ne i}f(\psi_j|1\le j\le N)=m_i(\psi_i)\;,
\vspace{2pt}\\\displaystyle
\rien\hspace{40pt}
\mbox{and}\quad\forall\psi\in<\Theta>_\Gamma^N\,,\;
\bigwedge_{i=1}^N\psi_i=\bot
\Longrightarrow f(\psi_i|1\le i\le N)=0\,.
\end{array}
%\right.
\end{array}
\end{equation}
\subsubsection{Approximation of the rule}
The definition of our rule needs the maximization of the entropy of $f$ under various constraints.
An algorithm for solving this maximization is proposed in section~\ref{DSmTb2:Dmb:Sec:4}.
The problem is solved by means of a variational method.
By the way, it may be interesting to have a more direct computation of the rule.
In particular, better computations of the rule could be obtained by approximating the optimization problem.
\\[5pt]
As soon as a solution is feasible, there are many ways to approximate the rules.
The main point is to approximate the entropy $H(f)=-\sum_{\psi,\eta}f(\psi,\eta)\ln f(\psi,\eta)$ by a function $\widetilde{H}(f)$ such that $\widetilde{H}(f)\simeq H(f)$\,.
Then, the rule is just rewritten:
\begin{equation}\label{dmb:dsmtb2:ch2:rule:eq:8}
\begin{array}{@{}l@{}}\displaystyle
\forall\phi\in<\Theta>_\Gamma\,,\;m_1\oplus m_2(\phi)=\sum_{\psi\wedge\eta=\phi}f_o(\psi,\eta)\;,
\vspace{5pt}\\\mbox{where:}
\vspace{5pt}\\
\rien\hspace{30pt}
%\left\{
\begin{array}{@{\;}l@{}}\displaystyle
f_o\in\arg\max_f\widetilde{H}(f)\;,
\vspace{5pt}\\\mbox{under constraints: }
\vspace{4pt}\\\displaystyle
\rien\hspace{40pt}
f\ge0\;,\quad
\sum_{\psi}f(\psi,\eta)=m_2(\eta)\;,\quad
\sum_{\eta}f(\psi,\eta)=m_1(\psi)\;,
\vspace{2pt}\\\displaystyle
\rien\hspace{40pt}
\mbox{and}\quad\forall\psi,\eta\in<\Theta>_\Gamma\,,\;
\psi\wedge\eta=\bot
\Longrightarrow f(\psi,\eta)=0\,.
\end{array}
%\right.
\end{array}
\end{equation}
An example of approximation is $\widetilde{H}(f)=-\sum_{\psi,\eta}f^2(\psi,\eta)$, which is obtained by a first order derivation of $\ln$.
Approximated rules will not be investigated in the present chapter.
\section{Implementation and properties}
\label{DSmTb2:Dmb:Sec:4}
This section is devoted to the development of basic properties of the rule EMR and to practical implementation on examples.
\subsection{Properties}
\label{DSmTb2:Dmb:Sec:4:subsec:prop}
\paragraph{Commutativity.}
Let $m_1$ and $m_2$ be two bbas over $\phi\in<\Theta>_\Gamma$\,.
By definition~(\ref{dmb:dsmtb2:ch2:rule:eq:5}), the fused bba $m_1\oplus m_2$ exists if and only if $m_2\oplus m_1$ exists.
Then $m_1\oplus m_2=m_2\oplus m_1$\,.
\paragraph{Neutral element.}
Define the bba of total ignorance $\nu$ by $\nu(\top)=1$\,.
Let $m$ be a bba over $\phi\in<\Theta>_\Gamma$\,.
Then the fused bba $m\oplus\nu$ exists, and $m\oplus\nu=m$\,.
\begin{description}
\item[Proof.]
Since $\sum_\phi f_o(\phi,\psi)=\nu(\psi)$ and $f_o\ge0$\,, it is deduced $f_o(\phi,\psi)=0$ unless $\psi=\top$.
\\
Now, since $\sum_\psi f_o(\phi,\psi)=m(\phi)$, it is concluded:
$$
f_o(\phi,\top)=m(\phi)
\quad\mbox{and}\quad
f_o(\phi,\psi)=0\mbox{ for }\psi\ne\top\;.
$$
This function satisfies the hypotheses of~(\ref{dmb:dsmtb2:ch2:rule:eq:5}), thus implying the existence of $m\oplus\nu$\,.\\
Then the result $m\oplus\nu=m$ is immediate.
\item[$\Box\Box\Box$]\rien
\end{description}
\paragraph{Belief enhancement.}
Let be given two bbas $m_1$ and $m_2$, and assume that there exists a fused bba $m_1\oplus m_2$ computed by~(\ref{dmb:dsmtb2:ch2:rule:eq:5}).
Denote by $\mathrm{Bel}_1\oplus\mathrm{Bel}_2$ the belief function related to $m_1\oplus m_2$\,.
Then:
\begin{equation}\label{dmb:dsmtb2:ch2:Implement:eq:1}
\mathrm{Bel}_1\oplus\mathrm{Bel}_2(\phi)\ge\max\bigl\{\mathrm{Bel}_1(\phi),\mathrm{Bel}_2(\phi)\bigr\}
\quad\mbox{for any }\phi\in<\Theta>_\Gamma\;.
\end{equation}
\begin{description}
\item[Proof.]\rien\\
\emph{Proof of $\mathrm{Bel}_1\oplus\mathrm{Bel}_2(\phi)\ge\mathrm{Bel}_1(\phi)$\,.}\\
Let $f_o$ be a function satisfying to~(\ref{dmb:dsmtb2:ch2:rule:eq:5}).
\\
Then $\mathrm{Bel}_1\oplus\mathrm{Bel}_2(\phi)=\sum_{\psi\subset\phi} m_1\oplus m_2(\psi)= \sum_{\eta\wedge\xi\subset\phi}f_o(\eta,\xi)$\,.
\\
In particular, $\mathrm{Bel}_1\oplus\mathrm{Bel}_2(\phi)\ge\sum_{\psi\subset\phi}\sum_{\eta}f_o(\psi,\eta)$\,.
\\
At last, $\mathrm{Bel}_1\oplus\mathrm{Bel}_2(\phi)\ge\sum_{\psi\subset\phi}m_1(\psi)=\mathrm{Bel}_1(\phi)$\,.
\\[5pt]
\emph{Conclusion.}
It is similarly proved $\mathrm{Bel}_1\oplus\mathrm{Bel}_2(\phi)\ge\mathrm{Bel}_2(\phi)$ and then the final result.
\item[$\Box\Box\Box$]\rien
\end{description}
\emph{Corollary.}
Let be given two bbas $m_1$ and $m_2$, and assume that there exists a fused bba $m_1\oplus m_2$ computed by~(\ref{dmb:dsmtb2:ch2:rule:eq:5}).
Let $\phi_1,\dots,\phi_n\in<\Theta>_\Gamma$ be such that $\phi_i\wedge\phi_j=\bot$ for any $i\ne j$\,.
Then the property \emph{$\sum_{i=1}^n\max\bigl\{\mathrm{Bel}_1(\phi_i),\mathrm{Bel}_2(\phi_i)\bigr\}\le 1$ is necessarily true.}
\\[5pt]
This result is a direct consequence of the belief enhancement.
It could be used as a criterion for proving the non existence of the fusion rule.
\paragraph{Associativity.}
The computed rule~$\oplus$ is not associative.
\begin{description}
\item[Proof.]
Consider the bbas $m_1$, $m_2$ and $m_3$  defined on the Boolean algebra $\{\bot, a, \neg a, \top\}$ by\,:
$$\left\{\begin{array}{@{\,}l@{}}\displaystyle
m_1(a)=m_2(a)=0.5\quad\mbox{and}\quad m_1(\top)=m_2(\top)=0.5\,,
\\\displaystyle
m_3(\neg a)=0.5\quad\mbox{and}\quad m_3(\top)=0.5\,.
\end{array}\right.$$
We are now comparing the fusions $(m_1\oplus m_2)\oplus m_3$ and $m_1\oplus (m_2\oplus m_3)$\,.
\\[5pt]
\emph{Computing $(m_1\oplus m_2)\oplus m_3$.}\\
First it is noticed that there is no possible conflict between $m_1$ and $m_2$, so that $m_1\oplus m_2$ could be obtained by means of the usual conjunctive rule:
$$
m_1\oplus m_2(a)=0.5\times0.5+0.5\times0.5+0.5\times0.5=0.75
\quad\mbox{and}\quad
m_1\oplus m_2(\top)=0.5\times0.5\;.
$$
As a consequence:
$$
\max\bigl\{\mathrm{Bel}_1\oplus\mathrm{Bel}_2(a),\mathrm{Bel}_3(a)\bigr\}
+
\max\bigl\{\mathrm{Bel}_1\oplus\mathrm{Bel}_2(\neg a),\mathrm{Bel}_3(\neg a)\bigr\}
=
0.75+0.5>1\;.
$$
It is concluded that \underline{$(m_1\oplus m_2)\oplus m_3$ does not exist.}
\\[5pt]
\emph{Computing $m_1\oplus (m_2\oplus m_3)$.}\\
It is known that $\mathrm{Bel}_2\oplus\mathrm{Bel}_3\ge\max\bigl\{\mathrm{Bel}_2,\mathrm{Bel}_3\bigr\}$ when $m_2\oplus m_3$ exists.
\\
Since $\max\bigl\{\mathrm{Bel}_2(a),\mathrm{Bel}_3(a)\bigr\}=\max\bigl\{\mathrm{Bel}_2(\neg a),\mathrm{Bel}_3(\neg a)\bigr\}=0.5$\,, it is deduced that necessarily $m_2\oplus m_3(a)=m_2\oplus m_3(\neg a)=0.5$
\\
It appears that $m_2\oplus m_3(a)=m_2\oplus m_3(\neg a)=0.5$ is actually a valid solution, since it is derived from $f_o$ such that $f_o(a,\top)=f_o(\top,\neg a)=0.5$ (zeroed on the other cases).
\\
It is also deduced by a similar reasoning that $m_1\oplus (m_2\oplus m_3)$ exists and is necessary defined by \underline{$m_1\oplus (m_2\oplus m_3)(a)=m_1\oplus (m_2\oplus m_3)(\neg a)=0.5$\,.}
\\[5pt]
The associativity thus fails clearly on this example.
\item[$\Box\Box\Box$]\rien
\end{description}
\paragraph{Compatibility with the probabilistic bound hypothesis.}
A temptation in evidence theories is to link the notion of probability with the notion of belief by means of the relation:
\begin{equation}\label{dmb:dsmtb2:ch2:Implement:eq:1:1}
\mathrm{Bel}(\phi)\le p(\phi)\quad\mbox{for any }\phi\in<\Theta>_\Gamma\;.
\end{equation}
In general, this relation is not compatible with the fusion rules.
\begin{quote}
For example, let us test Dempster-Shafer rule on the relation~(\ref{dmb:dsmtb2:ch2:Implement:eq:1:1})\\
Let be given $m_1$ and $m_2$ defined on $\{\bot,a,\neg a,\top\}$ by $m_1(a)=m_1(\neg a)=0.5$ and $m_2(a)=m_2(\top)=0.5$\,.
\\
It is deduced $\mathrm{Bel}_1(a)=\mathrm{Bel}_1(\neg a)=0.5$\,, $\mathrm{Bel}_2(a)=0.5$ and $\mathrm{Bel}_2(\neg a)=0$\,.
\\
The choice of $m_1$ and $m_2$ is thus compatible with the bound hypothesis~(\ref{dmb:dsmtb2:ch2:Implement:eq:1:1}), and it follows $p(a)=p(\neg a)=0.5$\,.
\\
Now, Dempster-Shafer rule implies $m_1\oplus m_2(a)=2/3$ and $m_1\oplus m_2(\neg a)=1/3$\,.
\\
Confronting $m_1\oplus m_2$ to~(\ref{dmb:dsmtb2:ch2:Implement:eq:1:1}), it comes $p(a)\ge 2/3$\,.
\\
This is contradictory with the previously obtained relation $p(\neg a)= 0.5$\,. 
\end{quote}
This difficulty is avoided by some theorists by saying that the notion of probability is dependent of the considered sensor, or that belief and probability are two separated notions.
\\[5pt]
In our opinion, probability should be considered as an absolute notion.
We will see in section~\ref{DSmTb2:Dmb:Sec:5}, that the belief could be considered as a probabilistic modal proposition.
Then there are two cases:
\begin{itemize}
\item If the modal propositions are not comparable to the propositions without modality, then there is no obvious relation between the belief and the probability.
This is particularly the case of the TBM paradigm.
\item If the modal propositions are comparable to the propositions without modality (axiom m.iii of section~\ref{DSmTb2:Dmb:Sec:5}), then the bound hypothesis~(\ref{dmb:dsmtb2:ch2:Implement:eq:1:1}) is recovered.
Moreover, the fusion rule EMR is then logically derived.
\end{itemize}
This anticipatory logical result has the following implication:
\begin{center}
\emph{The rule EMR is compatible with the bound hypothesis~(\ref{dmb:dsmtb2:ch2:Implement:eq:1:1}).}
\end{center}
But this result is already foreseeable from the property~(\ref{dmb:dsmtb2:ch2:Implement:eq:1}).
Indeed, property~(\ref{dmb:dsmtb2:ch2:Implement:eq:1}) makes impossible the construction of a counter-example like the previous one of this paragraph.
\paragraph{Idempotence.}
The rule is not idempotent, since it allways increases the precision of a bba.
However it will be idempotent, when the bba does not contain any imprecision (\emph{e.g.} a probability).
\\[5pt]
\emph{This obvious property is just illustrated on examples subsequently.}
\subsection{Algorithm}
\label{DSmTb2:Dmb:Sec:4:subsec:algo}
The optimization~(\ref{dmb:dsmtb2:ch2:rule:eq:5}) is convex and is not difficult.
The implemented algorithm is based on Rosen's gradient projection method.
Now, the gradient of $H(f)=\sum_{\psi,\eta}-f(\psi,\eta)\ln f(\psi,\eta)$ is characterized by:
$$
D_fH(f)=\sum_{\psi,\eta}-(1+\ln f(\psi,\eta))\,df(\psi,\eta)\,.
$$
Then, the algorithm follows the synopsis:
\begin{enumerate}
\item Find a feasible solution $f_0$ to the simplex:
$$
\begin{array}{@{}l@{}}\displaystyle
f_0\ge0\;,\quad
\sum_{\psi}f_0(\psi,\eta)=m_2(\eta)\;,\quad
\sum_{\eta}f_0(\psi,\eta)=m_1(\psi)\;,
\vspace{2pt}\\\displaystyle
\rien\hspace{40pt}
\mbox{and}\quad\forall\psi,\eta\in<\Theta>_\Gamma\,,\;
\psi\wedge\eta=\bot
\Longrightarrow f_0(\psi,\eta)=0\,.
\end{array}
$$
If such a solution does not exist, then stop: \emph{the fusion is not possible.}
\\[5pt]
Otherwise, set $t=0$ and continue on next step.
\item\label{DSmTb2:Dmb:Sec:4:subsec:algo:process:step1} Build $\Delta f_t$ by solving the linear program:
$$
\begin{array}{@{\;}l@{}}\displaystyle
\max_{\Delta f_t}\sum_{\psi,\eta}-(1+\ln f_t(\psi,\eta))\,\Delta f_t(\psi,\eta)\;,
\vspace{5pt}\\\mbox{under constraints: }
\vspace{4pt}\\\displaystyle
\rien\hspace{40pt}
f_t+\Delta f_t\ge0\;,\quad
\sum_{\psi}\Delta f_t(\psi,\eta)=\sum_{\eta}\Delta f_t(\psi,\eta)=0\;,
\vspace{2pt}\\\displaystyle
\rien\hspace{40pt}
\mbox{and}\quad\forall\psi,\eta\in<\Theta>_\Gamma\,,\;
\psi\wedge\eta=\bot
\Longrightarrow \Delta f_t(\psi,\eta)=0\,.
\end{array}
$$
\item Repeat $\Delta f_t:=\Delta f_t/2$
until $H(f_t+\Delta f_t)>H(f_t)$\,. 
\item Set $f_{t+1}=f_t+\Delta f_t$.
Then set $t:=t+1$\,.
\item Reiterate from step~\ref{DSmTb2:Dmb:Sec:4:subsec:algo:process:step1} until full convergence.
\end{enumerate}
The linear programing library Coin-LP has been used in our implementation.
\subsection{Example}
\label{DSmTb2:Dmb:Sec:4:subsec:example}
In this section is studied the fusion of bbas $m_i$ defined over $\mathcal{P}\Bigl(\{a,b,c\}\Bigr)$ by:
$$\left\{\begin{array}{@{\,}l@{}}
m_1(a)=\alpha_1\,,\ 
m_1(b)=0\,,\ 
m_1(c)=\gamma_1
\mbox{ and }
m_1(\{a,b,c\})=1-\alpha_1-\gamma_1\;,
\vspace{4pt}\\
m_2(a)=0\,,\ 
m_2(b)=\beta_2\,,\ 
m_2(c)=\gamma_2
\mbox{ and }
m_2(\{a,b,c\})=1-\beta_2-\gamma_2\;.
\end{array}\right.$$
This is a slight generalization of Zadeh's example.
The fusion $m_1\oplus m_2$ is solved by the algorithm, but also mathematically.
The solutions were identical by the both methods.
The results of fusion are presented for particular choices of the parameters $\alpha,\beta,\gamma$.
\paragraph{Mathematical solution.}
Assume that $f$ is is a function over $\mathcal{P}\Bigl(\{a,b,c\}\Bigr)^2$ verifying the conditions~(\ref{dmb:dsmtb2:ch2:rule:eq:6}) of the rule.
The marginal constraints say:
\begin{equation}\label{DSmTb2:Dmb:Sec:4:subsec:algo:eq0}
\left\{\begin{array}{@{\,}l@{}}\displaystyle
\sum_{B\subset\{a,b,c\}}f(a,B)=\alpha_1\;,
\ %\
\sum_{B\subset\{a,b,c\}}f(c,B)=\gamma_1\;,
\ %\
\sum_{B\subset\{a,b,c\}}f(\{a,b,c\},B)=1-\alpha_1-\gamma_1\;,
\vspace{7pt}\\\displaystyle
\sum_{A\subset\{a,b,c\}}f(A,b)=\beta_2\;,
\ %\
\sum_{A\subset\{a,b,c\}}f(A,c)=\gamma_2\;,
\ %\
\sum_{A\subset\{a,b,c\}}f(A,\{a,b,c\})=1-\beta_2-\gamma_2\;,
\vspace{7pt}\\\displaystyle
\sum_{B\subset\{a,b,c\}}f(A,B)=0\quad\mbox{and}\quad\sum_{A\subset\{a,b,c\}}f(A,B)=0\mbox{ in any other cases.}
\end{array}\right.
\end{equation}
Since $f(A,B)=0$ for any $(A,B)$ such that $A\cap B=\emptyset$ and $f\ge 0$, it is deduced that :
$$\mbox{\small$
f\Bigl(a,\{a,b,c\}\Bigr)\,,\ 
f\Bigl(\{a,b,c\},b\Bigr)\,,\ 
f(c,c)\,,\ 
f\Bigl(c,\{a,b,c\}\Bigr)\,,\ 
f\Bigl(\{a,b,c\},c\Bigr)
\mbox{ and }
f\Bigl(\{a,b,c\},\{a,b,c\}\Bigr)
$}$$
are the only values of $f$\,, which are possibly non zero.
Then the system~(\ref{DSmTb2:Dmb:Sec:4:subsec:algo:eq0}) reduces to the linear solution:
\begin{equation}\label{DSmTb2:Dmb:Sec:4:subsec:algo:eq1}
\left\{\begin{array}{@{\,}l@{}}
f\Bigl(a,\{a,b,c\}\Bigr)=\alpha_1\,,\ 
f\Bigl(\{a,b,c\},b\Bigr)=\beta_2\,,\ 
f(c,c)=\theta\,,\ 
f\Bigl(c,\{a,b,c\}\Bigr)=\gamma_1-\theta\,,
\vspace{3pt}\\
f\Bigl(\{a,b,c\},c\Bigr)=\gamma_2-\theta\,,\ 
f\Bigl(\{a,b,c\},\{a,b,c\}\Bigr)=1-\alpha_1-\beta_2-\gamma_1-\gamma_2+\theta\;,
\vspace{3pt}\\
f(A,B)=0\mbox{ for any other case.}
\end{array}\right.
\end{equation}
This solution depends on the only unknow parameter~$\theta$.
The optimal parameter $\theta_o$ is obtained by solving:
$$
\max_\theta\left(
h(\theta)+h(\gamma_1-\theta)+h(\gamma_2-\theta)+h(1-\alpha_1-\beta_2-\gamma_1-\gamma_2+\theta)
\right)
\quad
\mbox{where }
h(\tau)=-\tau\ln\tau\,.
$$
The function $f_o$ is then computed by using $\theta_o$ in~(\ref{DSmTb2:Dmb:Sec:4:subsec:algo:eq1}).
And $m_1\oplus m_2$ is of course obtained by $m_1\oplus m_2(\phi)=\sum_{\psi\wedge\eta=\phi}f_o(\psi,\eta)$\,.
\\[5pt]
It is sometimes impossible to find a valid parameter $\theta$.
The condition of existence is easily derived:
\begin{equation}\label{DSmTb2:Dmb:Sec:4:subsec:algo:eq2}
\theta\mbox{ exists}\quad\mbox{if and only if}\quad
\max\{0,\alpha_1+\beta_2+\gamma_1+\gamma_2-1\} \le \min\{\gamma_1,\gamma_2\}\;.
\end{equation}
When this condition is fulfilled, the optimal parameter is given by\footnote{Other cases are not compliant with the existence condition~(\ref{DSmTb2:Dmb:Sec:4:subsec:algo:eq2}).}:
\begin{equation}\label{DSmTb2:Dmb:Sec:4:subsec:algo:eq3}
\theta_o=\frac{\gamma_1\gamma_2}{1-\alpha_1-\beta_2}
\mbox{ when }\alpha_1+\beta_2<1\;,
\quad
\theta_o=0\mbox{ when }\alpha_1+\beta_2=1\;.
\end{equation}
Then, it is not difficult to check that $\theta_o$ is bounded accordingly to the existence condition:
$$
\max\{0,\alpha_1+\beta_2+\gamma_1+\gamma_2-1\}\le\theta_o\le\min\{\gamma_1,\gamma_2\} \;.
$$
\paragraph{Experimentation.}\rien\\
\emph{Zadeh's example.}\\
Zadeh's example is defined by $\alpha_1=\beta_2=0,99$ and $\gamma_1=\gamma_2=0.01$\,.
This fusion is \emph{rejected by EMR.}
\\[5pt]
More generally, assume $\gamma_1=1-\alpha_1$ and $\gamma_2=1-\alpha_2$\,.
Then the condition:
$$\max\{0,\alpha_1+\beta_2+\gamma_1+\gamma_2-1\} \le \min\{\gamma_1,\gamma_2\}$$
fails unless when $\gamma_1=\gamma_2=1$\,.
The case $\gamma_1=\gamma_2=1$ is trivial, since it means a perfect agreement between the two sources.
Thus, Zadeh's example is rejected by EMR, even if there is a negligible conflict between the two sources.
\\[5pt]
By the way, this is not surprising.
In Zadeh's example, the bbas are put on the singletons only.
Then, the bbas work like probabilities, thus defining an uncertainty but without any imprecision.
Since the information provided by the sources are free from any imprecision, there are only two possible cases:  either the information are the same, either some information are false.
\\[5pt]
Now, imagine again that our information come with a negligible conflict:
$$
m_1(a)=m_2(b)=\epsilon\quad\mbox{and}\quad m_1(c)=m_2(c)=1-\epsilon\;.
$$
This could indeed happen, when our information have been obtained from slightly distorted sources.
Now, it has been seen that EMR rejects this fusion.
Thus, we have to be cautious when using EMR and the following recommendation has to be considered:
\begin{quote}
If the sources of information are distorted, even slightly, these distortions {\bf have to} be encoded in the bbas by an additional imprecision.
\end{quote}
Typically, by weakening the bbas as follows:
$$\left\{\begin{array}{@{\,}l@{}}
m_1(a)=m_2(b)=\rho\epsilon\;,\ m_1(c)=m_2(c)=\rho(1-\epsilon)\;,
\vspace{4pt}\\
m_1\bigl(\{a,b,c\}\bigr)=m_2\bigl(\{a,b,c\}\bigr)=1-\rho
\;,\mbox{ with }\rho\le\frac1{1+\epsilon}\;,
\end{array}\right.$$
the fusion is again possible by means of EMR.
\\[5pt]
\emph{Extended example.}\\
The following table summarizes the fusion by EMR for various values of $\alpha,\beta,\gamma$\,:
$$
\begin{array}{@{}c|c||c|c||l@{}}
\alpha_1 & \gamma_1 & \beta_2 & \gamma_2 & m=m_1\oplus m_2
\\\hline
0.501&0&0.501&0& \mathrm{Rejection}
\\\hline
0.499&0&0.499&0& m(a)=m(b)=0.499\,,\ m(\{a,b,c\})=0.002
\\\hline
 0.3&0.1
&0.3&0.1
& m(a)=m(b)=0.3\,,\ m(c)=0.175\,,\ m(\{a,b,c\})=0.225
\\\hline
 0.3&0.05
&0.3&0.05
& m(a)=m(b)=0.3\,,\ m(c)=0.09375\,,\ m(\{a,b,c\})=0.30625
\\\hline
 0.3&0.01
&0.3&0.01
& m(a)=m(b)=0.3\,,\ m(c)=0.01975\,,\ m(\{a,b,c\})=0.38025
\end{array}$$
\emph{Comparison.}\\
In this last example, we compare the fusion by EMR and by Dempster-Shafer of the bbas $m_1$ and $m_2$ defined by:
$$\left\{\begin{array}{@{\,}l@{}}
m_1(a)=
m_1(\{a,b\})=
m_1(\{a,c\})=
m_1(\{b,c\})=
m_1(\{a,b,c\})=0.2\;,
\vspace{4pt}\\
m_2(a)=
m_2(\{a,b\})=
m_2(\{a,c\})=
m_2(\{b,c\})=
m_2(\{a,b,c\})=0.2\;.
\end{array}\right.$$
The following table presents the fusion results for Dempster-Shafer (DST) and for EMR:
$$\begin{array}{@{}c||c|c|c|c|c|c|c@{}}
A
&a&b&c
&\{a,b\}&\{a,c\}&\{b,c\}&\{a,b,c\}
\\\hline\hline
m_1\oplus m_2(A)\;/\;DST
& 0.390 & 0.087 & 0.087 & 0.131 & 0.131 & 0.131 & 0.043
\\\hline
m_1\oplus m_2(A)\;/\;EMR
& 0.411 & 0.093 & 0.093 & 0.107 & 0.107 & 0.153 & 0.036
\end{array}$$
We will not argue about the advantage of EMR compared to DST on such simple example.
The important point is to notice how the bba concentration subtly changes from DST to EMR.
In general, the belief enforcement of the small propositions is stronger in EMR.
But the case of proposition $\{b,c\}$ is different, since it is made weaker by DST than by EMR.
This is perhaps a consequence of a greater belief attraction of proposition $a$ compared to $b$ and $c$, during the fusion process.
\section{Logical ground of the rule}
\label{DSmTb2:Dmb:Sec:5}
This section justifies logically the definition~(\ref{dmb:dsmtb2:ch2:rule:eq:5}) of EMR.
This logical development, based on modal logics, is quite similar to what have been previously done in the DSmT book~1, chapter~8~\cite{DSmTb1:dmb}.
Actually, the modal operators will be applied to non modal proposition only for the simplicity of the exposition (the definitions of~\cite{DSmTb1:dmb} were more general), but there is no significant change in the logic.
Now, the reader could also refer to~\cite{DSmTb1:dmb} since it introduces the logic on examples.
Subsequently, the logic behind EMR will be exposed directly.
\\[5pt]
In~\cite{DSmTb1:dmb} the definition of the logic was axiomatic.
Since logic is not usual for many readers, such axiomatic definition is avoided in this presentation.
A model viewpoint is considered now.
More precisely, the description of our modal logic is made in the framework of a Boolean algebra.
Typically, a logical relation $\vdash \phi$, \emph{i.e. $\phi$ is proved}, is replaced by $\phi=\top$ in the Boolean algebra.
In the same way, a relation $\vdash \phi\rightarrow \psi$ is replaced by $\phi\subset\psi$.
Moreover, the modal relations are directly used within the Boolean algebra, \emph{i.e.} $\Box\phi\subset\phi$ is the Boolean counterpart of the logical axiom $\vdash\Box\phi\rightarrow\phi$.
We will not justify here the soundness of this Boolean framework, in regards to the modal propositions.
But such framework is of course tightly related to an implied Kripke's model.
By the way, it is also assumed that our model is complete for the logic.
But all these technical considerations should be ignored by most readers.
\paragraph{Notation.}
Subsequently, the proposition $\phi\setminus\psi$ is a shortcut for $\phi\wedge\neg\psi$\,.
Moreover, the notation $\phi\subsetneq\psi$ is used as a shortcut for $\phi\subset\psi$ and $\phi\ne\psi$ .
\subsection{The belief as a probability of a modal proposition}
\label{DSmTb2:Dmb:Sec:5:subsec:intro}
Many evidence theorists are particularly cautious, when comparing belief and probabilities.
On the first hand, there is a historical reason.
As new theories of uncertainty, the evidence theories had to confront the already existing theory of probability.
On the second hand, the sub-additivity property of the belief is often proposed as a counter-argument against the probabilistic hypothesis.
In this introductory subsection, it is seen that this property is easily fulfilled by means of modal propositions.
\paragraph{Sub-additivity, modality and probability.}
Subsequently, the belief $\mathrm{Bel}_i(\phi)$ of a proposition $\phi$ according to a sensor $i$ is interpreted as the probability $p(\Box_i\phi)$ of the modal proposition $\Box_i\phi$\,.
It is not the purpose of this paragraph to make a full justification of this interpretation, but rather to explain how the modal interpretation is able to overcome the sub-additivity property.
\\[5pt]
First at all, it is important to understand the meaning of the modal operator $\Box_i$\,.
The modal proposition $\Box_i\phi$ characterizes all the information, which sensor $i$ can associate for sure to the proposition $\phi$.
\begin{quote}
For example, the equation $\Box_i\phi=\top$ means: \emph{sensor $i$ considers $\phi$ as true in any configuration of the system}.
\end{quote}
Then, it is noticed that the modal propositions fulfill a \emph{logical sub-additivity} property:
\begin{equation}
\label{DSmTb2:Dmb:Sec:5:subsec:intro:eq:1}
(\Box_i\phi\vee\Box_i\psi)\subset\Box_i(\phi\vee\psi)\;.
\end{equation}
The converse is false in general.
This well-known property will be proved in the subsequent section.
It has a concrete interpretation.
The proposition $\Box_i\phi$ describes logically the information about $\phi$ which is granted as sure by the sensor $i$.
But there are information which sensor $i$ can attribute to $\phi\vee\psi$ for sure, but which cannot be attributed without ambiguity to $\phi$ or to $\psi$ alone.
These ambiguous information are exactly described by the non empty proposition $\Box_i(\phi\vee\psi)\setminus(\Box_i\phi\vee\Box_i\psi)$\,.
\\[5pt]
The important point now is that the logical sub-additivity directly implies the sub-additivity of the belief.
From~(\ref{DSmTb2:Dmb:Sec:5:subsec:intro:eq:1}), it is derived:
$$
p(\Box_i\phi\vee\Box_i\psi)\le p\bigl(\Box_i(\phi\vee\psi)\bigr)\;.
$$
Assume now that $\phi$ and $\psi$ are disjoint, \emph{i.e.} $\phi\wedge\psi=\bot$.
It is also hypothesized that $\Box_i$ is coherent, which implies $\Box_i\phi\wedge\Box_i\psi=\bot$ \emph{(coherent sensors will consider disjoint propositions as disjoint)}.
Then it comes:
$$
p(\Box_i\phi)+p(\Box_i\psi)\le p\bigl(\Box_i(\phi\vee\psi)\bigr)\;,
$$
and it is finally deduced:
$$
\mathrm{Bel}_i(\phi)+\mathrm{Bel}_i(\psi)\le\mathrm{Bel}_i(\phi\vee\psi)\;,
$$
from the modal definition of $\mathrm{Bel}_i$\,.
At last, we have recovered the sub-additivity of $\mathrm{Bel}_i$ from the logical sub-additivity of the modal operator $\Box_i$\,.
\\[5pt]
It appears that the sub-additivity is not incompatible with a probabilistic interpretation of the belief.
The sub-additivity seems rather related to a modal nature of the belief.
At the end of this paragraph, a last word have to be said about the TBM paradigm.
Beliefs as defined in~(\ref{f2k5:DSmTcont:Eq:1}), that is including the mass assignment of the empty set, is not sub-additive in TBM.
Only the truncated belief $\mathrm{Bel}_{S}$ is sub-additive.
This is a consequence of the possible non-zero mass assignment of the empty set.
By the way, there is also a modal interpretation of this fact.
It is seen subsequently that the TBM paradigm is equivalent to remove the coherence hypothesis about the modality $\Box_i$.
But the incoherence of $\Box_i$ allows $\Box_i\phi\wedge\Box_i\psi\ne\bot$ even when $\phi\wedge\psi=\bot$.
As a consequence, the sub-additivity cannot be recovered from the modal hypothesis either, when considering incoherent sensors.
\\[5pt]
This introduction has revealed some similar characteristics of the belief functions and the modal operators.
The forthcoming sections establish a complete correspondence between the evidence theories and the probabilistic modal logic.
%\paragraph{How the computed rule is compatible with a probabilistic hypothesis?}
%\ref{DSmTb2:Dmb:Sec:4:sub:BelEnh}
%
\subsection{Definition of the logic}
\label{DSmTb2:Dmb:Sec:5:subsec:def}
Let $\Theta$ be a set of atomic propositions.
Let $\Gamma$ be a set of Boolean constraints built from $\Theta$ and the Boolean operators.
The constraints of $\Gamma$ are describing the logical priors which are known about the atomic propositions.
Then, define $\mathcal{B}_\Gamma(\Theta)$ the Boolean algebra generated by $\Theta$ and compliant with the constraints of $\Gamma$.\footnote{$\mathcal{B}_\Gamma(\Theta)$ could be obtained by propagating the constraints of $\Gamma$ within the free Boolean algebra $\mathcal{B}(\Theta)$.}
The algebra $\mathcal{B}_\Gamma(\Theta)$ will describe the real world.
\\[5pt]
The real world will be observed by sensors.
Let $I$ be the set of the sensors.
The sensors of $I$ may be combined by pairs, thus constituting composite sensors.
The set of composite sensors, denoted $J$, is defined recursively as follows:
$$
I\subset J
\quad\mbox{and} \quad
(i,j)\in J\mbox{ for any }i,j\in J\;.
$$
Among the (composite) sensors of $J$, it will be assumed that some pairs of sensors are mutually independent.
\begin{equation}
\label{DSmTb2:Dmb:Sec:5:subsec:def:eq:1}
\mbox{\emph{For $i,j\in J$, the notation $i\times j$ indicates that the sensors $i$ and $j$ are independent.}}
\end{equation}
To any sensor $i\in J$ is associated a modal operators $\Box_i$ (for short, the notation $\Box_{i,j}$ will be used instead of $\Box_{(i,j)}$).
The modal operators $\Box_i$ will describe logically the information received and interpreted by sensor $i$.
The modal operators will act on any proposition of $\mathcal{B}_\Gamma(\Theta)$.
For any $\phi\in\mathcal{B}_\Gamma(\Theta)$\,, the proposition $\Box_i\phi$ is an interpretation of the real proposition $\phi$ by the sensor $i$.
It is noticed that $\Box_i\phi$ is not necessarily a proposition of $\mathcal{B}_\Gamma(\Theta)$, so that $\Box_i$ should be considered as an external operator.
In order to manipulate these modal propositions, we will consider the Boolean algebra $\mathcal{B}_\Gamma(\Theta,\Box)$ generated by $\Theta$ and $\Box_i\phi$ where $\phi\in\mathcal{B}_\Gamma(\Theta)$ and $i\in J$\,.
It is also assumed that the algebra $\mathcal{B}_\Gamma(\Theta,\Box)$ is compliant with the following constraints on the modal propositions:
\begin{description}
\item[m.i.] $\Box_i\top=\top$\,,
\item[m.ii.] $\Box_i(\neg\phi\vee\psi)\subset(\neg\Box_i\phi\vee\Box_i\psi)$\,, for any $\phi,\psi\in\mathcal{B}_\Gamma(\Theta)$ and $i\in I$\,,
\item[m.iii.] (optional) $\Box_i\phi\subset\phi$\,, for any $\phi\in\mathcal{B}_\Gamma(\Theta)$ and $i\in I$\,,
\item[m.iv.] $\Box_{i}\phi\subset\Box_{i,j}\phi$\,, for any $\phi\in\mathcal{B}_\Gamma(\Theta)$ and $i,j\in I$\,,
\item[m.indep.] $\Box_{i,j}\phi\subset(\Box_{i}\phi\vee\Box_{j}\phi)$\,, for any $\phi\in\mathcal{B}_\Gamma(\Theta)$ and $i,j\in I$ such that $i\times j$\,.
\end{description}
Together with the axioms $m.\ast$, the algebra $\mathcal{B}_\Gamma(\Theta,\Box)$ is a model characterizing our modal logic.
It is necessary to explain the signification of these axioms.
\\[3pt]
Axiom m.i explains that the sensors hold any tautology (trivially true propositions) for true.
\\[3pt]
Axiom m.ii is the model counterpart of the axiom $\vdash\Box_i(\phi\rightarrow\psi)\rightarrow(\Box_i\phi\rightarrow\Box_i\psi)$ of modal logic, which is a \emph{modus ponens} encoded within the modality.
In other word, axiom m.ii just says that the sensors make logical deductions.
Together, axioms m.i and m.ii express that the sensors are reasoning logically. 
\\[3pt]
Axiom m.iii says that the sensors always say the truth.
More precisely, it says that \emph{$\phi$ is true} whenever \emph{sensor $i$ considers $\phi$ as true}.
This axiom implies the coherence of the sensors.
By the way, it is probably stronger than the coherence hypothesis.
Axiom m.iii is considered as optional.
It will be suppressed in the case of the TBM paradigm, but used otherwise.
\\[3pt]
Axiom m.iv says that the knowledge increases with the number of sensors used.
More precisely, m.iv means that the surety of a proposition $\phi$ is greater according to $i,j$ than according to $i$ only.
Although this axiom seems quite natural, it is noticed that it is not necessarily compatible with fusion rules involving a redistribution of the conflict.
\\[3pt]
Axiom m.indep expresses a logical independence of the sensors.
Assuming $i,j$ to be independent sensors (\emph{i.e.} $i\times j$), the axiom m.indep then says that the information obtained from the \emph{joint} sensor $(i,j)$ could be obtained separately from one of the sensors $i$ or $j$.
In other word, there is no possible interaction between the sensors $i$ and $j$ during the observation process.
\\[5pt]
m.i, m.ii and m.iii are typically the axioms of the system T of modal logic.
Before concluding this subsection, some useful logical theorems are now derived.
\subsubsection{Useful theorems}
\label{DSmTb2:Dmb:Sec:5:subsec:useTh}
The following theorems will be deduced without the help of the optional axiom m.iii.
The axioms used for each theorem will be indicated at the end of the proof.
\paragraph{The modality is non decreasing.} Let $i\in J$ and $\phi,\psi\in\mathcal{B}_\Gamma(\Theta)$\,.
Then:
$$
\phi\subset\psi
\quad\mbox{implies}\quad
\Box_i\phi\subset\Box_i\psi
\;.
$$
\begin{description}
\item[Proof.]
$\phi\subset\psi$ implies $\neg\phi\vee\psi=\top$\,.
\\
By applying axiom m.i, it comes $\Box_i\bigl(\neg\phi\vee\psi\bigr)=\top$\,.
\\
Now, m.ii implies $\Box_i\bigl(\neg\phi\vee\psi\bigr)\subset \neg\Box_i\phi\vee\Box_i\psi$\,.
\\
Consequently $\neg\Box_i\phi\vee\Box_i\psi=\top$\,.
\\
As a conclusion, $\Box_i\phi\subset\Box_i\psi$\,.
\item[$\Box\Box\Box$]\rien
\\ The proof requested the axioms m.i and m.ii.
\end{description}
\paragraph{Modality and conjunction.} $(\Box_i\phi\wedge\Box_i\psi)=\Box_i(\phi\wedge\psi)$ for any $\phi,\psi\in\mathcal{B}_\Gamma(\Theta)$ and $i\in J$\,.
\begin{description}
\item[Proof.]\rien\\
\emph{Proof of }$(\Box_i\phi\wedge\Box_i\psi)\subset\Box_i(\phi\wedge\psi)$\,.
\\
Since $\neg\phi\vee\neg\psi\vee(\phi\wedge\psi)=\top$\,, axiom m.i implies $\Box_i\bigl(\neg\phi\vee\neg\psi\vee(\phi\wedge\psi)\bigr)=\top$\,.
\\
Now, m.ii implies $\Box_i\bigl(\neg\phi\vee\neg\psi\vee(\phi\wedge\psi)\bigr)\subset \neg\Box_i\phi\vee\neg\Box_i\psi\vee\Box_i(\phi\wedge\psi)$\,.
\\
Consequently $\neg\Box_i\phi\vee\neg\Box_i\psi\vee\Box_i(\phi\wedge\psi)=\top$\,.
\\
As a conclusion, $(\Box_i\phi\wedge\Box_i\psi)\subset\Box_i(\phi\wedge\psi)$\,.
\\[5pt]
\emph{Proof of }$(\Box_i\phi\wedge\Box_i\psi)\supset\Box_i(\phi\wedge\psi)$\,.
\\
Since $\Box_i$ is non decreasing, it is proved $\Box_i(\phi\wedge\psi)\subset\Box_i\phi$ and $\Box_i(\phi\wedge\psi)\subset\Box_i\psi$.
\item[$\Box\Box\Box$]\rien
\\ The proof requested the axioms m.i and m.ii.
\end{description}
\paragraph{Modality and disjunction.}
$(\Box_i\phi\vee\Box_i\psi)\subset\Box_i(\phi\vee\psi)$ for any $\phi,\psi\in\mathcal{B}_\Gamma(\Theta)$ and $i\in J$\,.
In general, $(\Box_i\phi\vee\Box_i\psi)\ne\Box_i(\phi\vee\psi)$\,.
\begin{description}
\item[Proof.]\rien\\
\emph{Proof of }$(\Box_i\phi\vee\Box_i\psi)\subset\Box_i(\phi\vee\psi)$\,.
\\
Since $\Box_i$ is non decreasing, it is proved $\Box_i\phi\subset\Box_i(\phi\vee\psi)$ and $\Box_i\psi\subset\Box_i(\phi\vee\psi)$\,.
\\
Then the result.
\\[5pt] 
\emph{A counter-example for }$(\Box_i\phi\vee\Box_i\psi)=\Box_i(\phi\vee\psi)$\emph{ needs the construction of a dedicated model of the logic.
This construction is outside the scope of the chapter.
Readers interested in models constructions for modal logics could refer to \cite{blackburn}.}
\item[$\Box\Box\Box$]\rien
\\ The proof requested the axioms m.i and m.ii.
\end{description}
\paragraph{Conjunction of heterogeneous modalities.}
$(\Box_i\phi\wedge\Box_j\psi)\subset\Box_{i,j}(\phi\wedge\psi)$ for any \mbox{$\phi,\psi\in\mathcal{B}_\Gamma(\Theta)$} and $i,j\in J$\,.
\emph{In other words, if sensor $i$ asserts $\phi$ and sensor $j$ asserts $\psi$, then the fused sensor asserts $\phi\wedge\psi$\,.}
\begin{description}
\item[Proof.]
Axioms m.iv says $\Box_i\phi\subset\Box_{i,j}\phi$ and $\Box_j\psi\subset\Box_{i,j}\psi$\,.\\
Since $(\Box_{i,j}\phi\wedge\Box_{i,j}\psi)=\Box_{i,j}(\phi\wedge\psi)$\,, it is deduced $(\Box_i\phi\wedge\Box_j\psi)\subset\Box_{i,j}(\phi\wedge\psi)$\,.
\item[$\Box\Box\Box$]\rien
\\ The proof requested the axioms m.i, m.ii and m.iv.
\end{description}
\paragraph{Disjunction of heterogeneous modalities.}
$(\Box_i\phi\vee\Box_j\phi)\subset\Box_{i,j}\phi$ for any \mbox{$\phi\in\mathcal{B}_\Gamma(\Theta)$} and $i,j\in J$\,.
\emph{In other words, if sensor $i$ or sensor $j$ assert $\phi$, then the fused sensor asserts $\phi$\,.}
\begin{description}
\item[Proof.]
Axioms m.iv says $\Box_i\phi\subset\Box_{i,j}\phi$ and $\Box_j\phi\subset\Box_{i,j}\phi$\,.\\
Then, the result is immediate.
\item[$\Box\Box\Box$]\rien
\\ The proof requested the axiom m.iv.
\end{description}
The converse of this property is obtained by means of axiom m.indep, when the sensors $i,j$ are independent.
\subsection{Fusion rule}
\label{DSmTb2:Dmb:Sec:5:subsec:rule}
The purpose of this subsection is to derive logically the fusion rule on $<\Theta>_\Gamma$, the pre-Boolean algebra generated by $\Theta$ within the Boolean algebra $\mathcal{B}_\Gamma(\Theta)$.\footnote{In this case $\Gamma$ may contain constraints outside $<\Theta>$. But this is the same notion of pre-Boolean algebra discussed earlier.}
In a first step, the fusion will be derived in a strict logical acceptation, by means of the modal operators.
In a second step, the notion of belief is also introduced by means of probabilistic modal propositions.
But as a preliminary, we are beginning by introducing the notion of partitions.
\subsubsection{Preliminary definitions.}
\label{DSmTb2:Dmb:Sec:5:subsec:rule:preldef}
\paragraph{Partition.}
Let $\Pi\subset\mathcal{B}(\Theta)$ be a set of propositions.
The set $\Pi$ is a partition of $\top$ if it satisfies the following properties:
\begin{itemize}
\item The propositions of $\Pi$ are exclusive (\emph{i.e.} disjoint): $\phi\wedge\psi=\bot$ for any $\phi,\psi\in \Pi$ such that $\phi\ne\psi$\,,
\item The propositions of $\Pi$ are exhaustive: $\bigvee_{\phi\in \Pi}\phi=\top$\,.
\end{itemize}
Notice that $\Pi$ may contain $\bot$, in this definition.
\paragraph{Partition and complement.}
Let $\Pi$ be a partition and let $A\subset\Pi$ and $B\subset\Pi$.
Then:
$$
\left(\bigvee_{\phi\in A}\phi\right)\setminus\left(\bigvee_{\phi\in B}\phi\right)
=\left(\bigvee_{\phi\in A\setminus B}\phi\right)\;.
$$
This property just tells that the Boolean algebra generated by $\Pi$ is isomorph to the Boolean structure implied by the set $\Pi$.
The proof of this result is obvious from the definition.
\paragraph{Partitions and probabilities.}
Partitions are useful since they make possible the definition of a probability by means of elementary density.
More precisely, for any partition $\Pi$ and any subset $A\subset \Pi$, the probability of the proposition $\bigvee_{\phi\in A}\phi$ is given by:
$$
p\left(\bigvee_{\phi\in A}\phi\right)=\sum_{\phi\in A}p(\phi)\;.
$$
This property will be particularly useful subsequently for linking the logical fusion to the belief fusion.
\paragraph{Joint partitions.}
Let $\Pi$ and $\Lambda$ be two partitions of $\top$.
Let $\Gamma=\{\phi\wedge\psi/\phi\in \Pi\mbox{ and }\psi\in \Lambda\}$ be the set of joint propositions obtained from $\Pi$ and $\Lambda$.
Then $\Gamma$ is a partition.
\begin{description}
\item[Proof.]
Let $\phi,\phi'\in \Pi$ and $\psi,\psi'\in \Lambda$\,, such that $(\phi,\psi)\ne(\phi',\psi')$.
\\
The exclusivity of $(\phi\wedge\psi)$ and $(\phi'\wedge\psi')$ is a direct consequence of:
$$
(\phi\wedge\psi)\wedge(\phi'\wedge\psi')=(\phi\wedge\phi')\wedge(\psi\wedge\psi')=\bot\;.
$$
The exhaustivity is derived from:
$$
\bigvee_{\phi\in \Pi}\;\;\bigvee_{\psi\in \Lambda}(\phi\wedge\psi)=\left(\bigvee_{\phi\in \Pi}\phi\right)\wedge \left(\bigvee_{\psi\in \Lambda}\psi\right)=\top\wedge\top=\top\;.
$$
\item[$\Box\Box\Box$]\rien
\end{description}
\emph{Corollary of the proof.}
$(\phi\wedge\psi)=(\phi'\wedge\psi')$ and $(\phi,\psi)\ne(\phi',\psi')$ imply $(\phi\wedge\psi)=(\phi'\wedge\psi')=\bot$\,.
\\[5pt]
This corollary will be useful for the computation of $\phi^{(i,j)}$\,, subsequently.
\subsubsection{Logical fusion}
\paragraph{Definition of the logical fusion.}
Logically, the information provided by the sensor \mbox{$i\in J$} is described by the modal propositions $\Box_i\phi$, where $\phi\in<\Theta>_\Gamma$.
The propositions of \mbox{$\mathcal{B}_\Gamma(\Theta)\setminus<\Theta>_\Gamma$} are not considered explicitly, since our \emph{discernment} is restricted to $<\Theta>_\Gamma$.
\\[5pt]
Let $i,j\in J$ be two sensors which are independent, \emph{i.e.} such that $i\times j$.
The fusion of $i$ and $j$ is simply defined as the composite sensor $(i,j)$.
Now arises the following issue: \emph{How to characterize the fused information \mbox{$\Box_{i,j}\phi$} from the primary information \mbox{$\Box_i\phi$} and \mbox{$\Box_j\phi$}\,?}
In order to solve this question, we introduce first the notion of \emph{basic propositional assignments} which constitute the elementary logical components of the information.
\paragraph{Definition of the basic propositional assignments.}
Let $i\in J$ be a sensor.
The basic propositional assignments (bpa) related to sensor $i$ are the modal propositions $\phi^{(i)}$ defined for any $\phi\in<\Theta>_\Gamma$ by:
\begin{equation}
\label{DSmTb2:Dmb:Sec:5:subsec:rule:eq:1}
\phi^{(i)}=\Box_i\phi\setminus\left(
\bigvee_{\psi\in<\Theta>_\Gamma:\psi\subsetneq\phi}\Box_i\psi
\right)\;.
\end{equation}
The bpa $\phi^{(i)}$ is the logical information, which sensor $i$ attributes to proposition $\phi$ intrinsically.
The information of $\phi^{(i)}$ cannot be attributed to smaller propositions than $\phi$.
\\[5pt]
Subsequently, the bpas appear as essential tools for characterizing the fusion rule.
\paragraph{Logical properties of the bpa.}
\rien\\
\emph{Exclusivity.}
The bpas $\phi^{(i)}$, where $\phi\in<\Theta>_\Gamma$, are exclusive for any given sensor $i\in J$\,:
\begin{equation}
\label{DSmTb2:Dmb:Sec:5:subsec:rule:eq:2}
\forall \phi,\psi\in<\Theta>_\Gamma\,,\;
\phi\ne\psi\Rightarrow \phi^{(i)}\wedge\psi^{(i)}=\bot\;.
\end{equation}
\begin{description}
\item[Proof.]
From the definition, it is deduced:
$$
\phi^{(i)}\wedge\psi^{(i)}=\Box_i(\phi\wedge\psi)\wedge\left(
\bigwedge_{\eta\in<\Theta>_\Gamma:\eta\subsetneq\phi}\neg\Box_i\eta
\right)\wedge\left(
\bigwedge_{\eta\in<\Theta>_\Gamma:\eta\subsetneq\psi}\neg\Box_i\eta
\right)\;.
$$
Since $\phi\wedge\psi\subsetneq\phi$ or $\phi\wedge\psi\subsetneq\psi$ when $\phi\neq\psi$, it comes $\phi^{(i)}\wedge\psi^{(i)}=\bot$\,.
\item[$\Box\Box\Box$]\rien
\end{description}
%\rien\\
\emph{Exhaustivity.}
The bpas $\phi^{(i)}$, where $\phi\in<\Theta>_\Gamma$, are exhaustive for any given sensor $i\in J$\,:
\begin{equation}
\label{DSmTb2:Dmb:Sec:5:subsec:rule:eq:3}
\bigvee_{\psi\in<\Theta>_\Gamma:\psi\subset\phi}
\psi^{(i)}=\Box_i\phi\;,\quad
\mbox{and in particular:}
\bigvee_{\psi\in<\Theta>_\Gamma}\psi^{(i)}=\top\;.
\end{equation}
\begin{description}
\item[Proof.]
The proof is recursive.
\\
It is first noticed that $\Box_i\bot=\bot^{(i)}$\,.
\\[5pt]
Now, let $\phi\in<\Theta>_\Gamma$ and assume $\bigvee_{\eta\subset\psi}
\eta^{(i)}=\Box_i\psi$ for any $\psi\subsetneq\phi$\,.
\\
Then:
$$
\bigvee_{\psi\subset\phi}\psi^{(i)}=
\phi^{(i)}\vee\left(\bigvee_{\psi\subsetneq\phi}\;\;\bigvee_{\eta\subset\psi}\eta^{(i)}\right)=
\phi^{(i)}\vee\left(\bigvee_{\psi\subsetneq\phi}\Box_i\psi\right)\,.
$$
It follows $\bigvee_{\psi\subset\phi}\psi^{(i)}=
\left(\Box_i\phi\setminus\left(
\bigvee_{\psi\subsetneq\phi}\Box_i\psi
\right)\right)\vee\left(\bigvee_{\psi\subsetneq\phi}\Box_i\psi\right)=\Box_i\phi\vee\left(\bigvee_{\psi\subsetneq\phi}\Box_i\psi\right)\,.$ 
\\
Since $\Box_i$ is non decreasing, it is deduced $\bigvee_{\psi\subset\phi}\psi^{(i)}=\Box_i\phi$\,.
\item[$\Box\Box\Box$]\rien
\end{description}
%\rien\\
%
\emph{Partition.}
Being both disjoint and exhaustive, the bpas $\phi^{(i)}$, where $\phi\in<\Theta>_\Gamma$, constitute a partition of $\top$\,.
\rien\\[7pt]
\emph{Joint partition.}
Let $i,j\in I$.
The propositions $\phi^{(i)}\wedge\psi^{(j)}$, where $\phi,\psi\in<\Theta>_\Gamma$, constitute a partition of $\top$\,.
\paragraph{Computing the fusion.}
Let $i,j\in J$ be such that $i\times j$\,.
Then, the following property holds for any $\phi\in<\Theta>_\Gamma$\,:
\begin{equation}\label{DSmTb2:Dmb:Sec:5:subsec:rule:eq:4}
\phi^{(i,j)}=\bigvee_{\psi,\eta\in<\Theta>_\Gamma:\psi\wedge\eta=\phi}
\left(\psi^{(i)}\wedge\eta^{(j)}\right)\;.
\end{equation}
\begin{description}
\item[Proof.]\rien\\
\emph{Lemma.}
$$
\Box_{i,j}\phi=\bigvee_{\psi\wedge\eta\subset\phi}(\Box_{i}\psi\wedge\Box_{j}\eta)
=\bigvee_{\psi\wedge\eta\subset\phi}\left(\psi^{(i)}\wedge\eta^{(j)}\right)\;.
$$
\begin{description}
\item[Proof of lemma.]
From the property $\Box_i\psi\wedge\Box_j\eta\subset\Box_{i,j}(\psi\wedge\eta)$ of section~\ref{DSmTb2:Dmb:Sec:5:subsec:useTh} and the non decreasing property of $\Box_{i,j}$, it is deduced:
$$
\bigvee_{\psi\wedge\eta\subset\phi}(\Box_i\psi\wedge\Box_j\eta)\subset\Box_{i,j}\phi\;.
$$
Now, the axiom m.indep implies $\Box_{i,j}\phi\subset(\Box_{i}\phi\vee\Box_{j}\phi)$ and then:
$$
\Box_{i,j}\phi\subset\bigl((\Box_{i}\phi\wedge\Box_{j}\top)\vee((\Box_{i}\top\wedge\Box_{j}\phi)\bigr)\;.
$$
As a consequence, $\bigvee_{\psi\wedge\eta\subset\phi}(\Box_i\psi\wedge\Box_j\eta)=\Box_{i,j}\phi$\,.
\\[5pt]
Now, since $\Box_i\psi=\bigvee_{\xi\subset\psi}\xi^{(i)}$ and $\Box_j\eta=\bigvee_{\zeta\subset\eta}\zeta^{(j)}$ (refer to the exhaustivity property), it comes also:
$$
\Box_{i,j}\phi=\bigvee_{\psi\wedge\eta\subset\phi}\;\;
\bigvee_{\xi\subset\psi}\;\;\bigvee_{\zeta\subset\eta} \left(\xi^{(i)}\wedge\zeta^{(j)}\right)
=\bigvee_{\xi\wedge\zeta\subset\phi}\;\;\left(\xi^{(i)}\wedge\zeta^{(j)}\right)\;.
$$
\item[$\Box\Box$]\rien
\end{description}
From the definition of the bpa, it is deduced:
$$
\phi^{(i,j)}=\Box_{i,j}\phi\setminus\left(
\bigvee_{\psi\subsetneq\phi}\Box_{i,j}\psi
\right)=\bigvee_{\eta\wedge\xi\subset\phi}\left(\eta^{(i)}\wedge\xi^{(j)}\right)
\setminus\left(
\bigvee_{\eta\wedge\xi\subsetneq\phi}\left(\eta^{(i)}\wedge\xi^{(j)}\right)
\right)\;.
$$
Now, since the propositions $\eta^{(i)}\wedge\xi^{(j)}$ constitute a partition (and taking into account the \emph{corollary of the proof} in section~\ref{DSmTb2:Dmb:Sec:5:subsec:rule:preldef}), it comes:
$$
\phi^{(i,j)}=\bigvee_{\eta\wedge\xi=\phi}\left(\eta^{(i)}\wedge\xi^{(j)}\right)
\;.
$$
\item[$\Box\Box\Box$]\rien
\end{description}
\paragraph{Conclusion.}
The sensors $i,j\in J$ being independent, the fused sensor $(i,j)$ is computed from $i$ and $j$ accordingly to the following process:
\begin{itemize}
\item Build $\phi^{(i)}=\Box_i\phi\setminus\left(
\bigvee_{\psi\in<\Theta>_\Gamma:\psi\subsetneq\phi}\Box_i\psi
\right)=\Box_i\phi\setminus\left(
\bigvee_{\psi\in<\Theta>_\Gamma:\psi\subsetneq\phi}\psi^{(i)}
\right)$ and
\\
$\phi^{(j)}=\Box_j\phi\setminus\left(
\bigvee_{\psi\in<\Theta>_\Gamma:\psi\subsetneq\phi}\Box_j\psi
\right)=\Box_j\phi\setminus\left(
\bigvee_{\psi\in<\Theta>_\Gamma:\psi\subsetneq\phi}\psi^{(j)}
\right)$ for any $\phi\in<\Theta>_\Gamma$\,,
\item Compute $\phi^{(i,j)}=\bigvee_{\eta,\xi\in<\Theta>_\Gamma: \eta\wedge\xi=\phi}\left(\eta^{(i)}\wedge\xi^{(j)}\right)$ for any $\phi\in<\Theta>_\Gamma$\,,
\item Derive $\Box_{i,j}\phi=\bigvee_{\psi\in<\Theta>_\Gamma:\psi\subset\phi}
\psi^{(i,j)}$ for any $\phi\in<\Theta>_\Gamma$\,.
\end{itemize}
Obviously, this process is almost identical to the computation of the fused belief $\mathrm{Bel}_i\oplus\mathrm{Bel}_j$ in free DSmT or in the TBM paradigm (while including the empty proposition in the definition of the belief function):
\begin{itemize}
\item Set $m_i(\phi)=\mathrm{Bel}_i(\phi)-\sum_{\psi\subsetneq\phi}m_i(\psi)$
and $m_j(\phi)=\mathrm{Bel}_j(\phi)-\sum_{\psi\subsetneq\phi}m_j(\psi)$\,,
\item Compute $m_i\oplus m_j(\phi)=\sum_{\eta\wedge\xi=\phi}m_i(\eta)m_j(\xi)$\,,
\item Get back $\mathrm{Bel}_{i,j}(\phi)=\sum_{\psi\subset\phi}m_i\oplus m_j(\psi)$\,.
\end{itemize}
It is yet foreseeable that the $m_i\oplus m_j(\phi)$ could be interpreted as $p(\Box_{i,j}\phi)$ owing to some additional hypotheses about the \emph{probabilistic independence} of the proposition.
This idea will be combined with the entropic maximization method described in section~\ref{DSmTb2:Dmb:Sec:3}, resulting in a logically interpreted fusion rule for the evidence theories.
\\[5pt]
For now, we are discussing about the signification of optional axiom m.iii which has not been used until now.
\paragraph{The consequence of axiom m.iii.}
Axiom m.iii says $\Box_i\phi\subset\phi$ and in particular implies $\Box_i\bot\subset\bot$ and then $\Box_i\bot=\bot$\,.
Thus, there are two important properties related to m.iii:
\begin{itemize}
\item It establishes a comparison of the propositions $\phi$ and their interpretation $\Box_i\phi$ by means of $\Box_i\phi\subset\phi$\,,
\item It makes the sensors $i$ coherent by implying $\Box_i\bot=\bot$\,.
\end{itemize}
By removing m.iii, the incoherence $\Box_i\bot\ne\bot$ is made possible, and this has a fundamental interpretation in term of evidence theories.
\begin{itemize}
\item Allowing the incoherence $\Box_i\bot\ne\bot$ is a logical counterpart of the TBM paradigm,
\item Hypothesizing the coherence $\Box_i\bot=\bot$ is a logical counterpart of the DSmT or DST paradigm.
\end{itemize}
Next section establishes the connection between the logical fusion and the belief fusion.
\subsubsection{From logical to belief fusion}
Subsequently, we are assuming that a probability $p$ is defined over the Boolean algebra $\mathcal{B}_\Gamma(\Theta,\Box)$.
This probability is known partially by means of the sensors.
For any $i\in J$ and any $\phi\in<\Theta>_\Gamma$ are then defined:
\begin{itemize}
\item The belief $\mathrm{Bel}_i(\phi)=p(\Box_i\phi)$\,,
\item The basic belief assignment $m_i(\phi)=p\bigl(\phi^{(i)}\bigr)$\,.
\end{itemize}
For any $i,j\in J$ such that $i\times j$ (independent sensors), the fused bba and belief are defined by: 
\begin{equation}
\label{DSmTb2:Dmb:Sec:5:subsec:logtobelfus:eq:1}
m_i\oplus m_j=m_{i,j}
\quad\mbox{and}\quad
\mathrm{Bel}_i\oplus\mathrm{Bel}_j=\mathrm{Bel}_{i,j}\;.
\end{equation}
The propositions $\phi^{(i)}$ constituting a partition of $\top$, the logical property
\begin{equation}
\label{DSmTb2:Dmb:Sec:5:subsec:logtobelfus:eq:2}
\phi^{(i)}=\Box_i\phi\setminus\left(
\bigvee_{\psi\in<\Theta>_\Gamma:\psi\subsetneq\phi}\psi^{(i)}
\right)
\end{equation}
implies:
$$
m_i(\phi)=\mathrm{Bel}_i(\phi)-\sum_{\psi\subsetneq\phi}m_i(\psi)\;.
$$
From the exhaustivity property, \emph{i.e.} $\Box_{i}\phi=\bigvee_{\psi\subset\phi}
\psi^{(i)}$, is derived:
$$
\mathrm{Bel}_{i}(\phi)=\sum_{\psi\subset\phi}m_i(\psi)\;.
$$
By the way, two fundamental properties of evidence theories have been recovered from our logical approach.
Now, the remaining question is about the fusion rule.
\\[5pt]
From the definition and the computation of $\phi^{(i,j)}$, it is deduced:
$$
m_i\oplus m_j(\phi)=p\bigl(\phi^{(i,j)}\bigr)=p\left(
\bigvee_{\eta\wedge\xi=\phi}\left(\eta^{(i)}\wedge\xi^{(j)}\right)
\right)\;.
$$
Since the propositions $\eta^{(i)}\wedge\xi^{(j)}$ are constituting a partition (and owing to the \emph{corollary of the proof} in section~\ref{DSmTb2:Dmb:Sec:5:subsec:rule:preldef}), it is obtained:
\begin{equation}
\label{DSmTb2:Dmb:Sec:5:subsec:logtobelfus:eq:3}
m_i\oplus m_j(\phi)=\sum_{\eta\wedge\xi=\phi}p\left(\eta^{(i)}\wedge\xi^{(j)}
\right)\;.
\end{equation}
It is not possible to reduce~(\ref{DSmTb2:Dmb:Sec:5:subsec:logtobelfus:eq:3}) anymore, without an additional hypothesis.
In order to compute $p\left(\eta^{(i)}\wedge\xi^{(j)}\right)$\,, the independence of sensors $i$ and $j$ will be again instrumental.
But this time, the independence is considered from an entropic viewpoint, and the probabilities $p\left(\eta^{(i)}\wedge\xi^{(j)}\right)$ are computed by maximizing the entropy of $p$ over the propositions  $\eta^{(i)}\wedge\xi^{(j)}$.
Denoting $\mathcal{P}(\Theta)=\mathcal{P}\bigl(\mathcal{B}_\Gamma(\Theta)\bigr)$ the set of all probabilities over $\mathcal{B}_\Gamma(\Theta)$, the probabilities $p\left(\eta^{(i)}\wedge\xi^{(j)}\right)$ are obtained by means of the program:
\begin{equation}
\label{DSmTb2:Dmb:Sec:5:subsec:logtobelfus:eq:4}
\begin{array}{@{}l@{}}\displaystyle
p\in\max_{q\in\mathcal{P}(\Theta)}\sum_{\eta,\xi\in<\Theta>_\Gamma}-q\left(\eta^{(i)}\wedge\xi^{(j)}\right) \ln q\left(\eta^{(i)}\wedge\xi^{(j)}\right)
\;,
\vspace{7pt}\\\displaystyle
\mbox{under constraints:}\quad
q\left(\phi^{(i)}\right)=m_i(\phi)
\quad\mbox{and}\quad
q\left(\phi^{(j)}\right)=m_j(\phi)
\quad\mbox{for any }\phi\in<\Theta>_\Gamma\;.
\end{array}
\end{equation}
Combining~(\ref{DSmTb2:Dmb:Sec:5:subsec:logtobelfus:eq:3}) and~(\ref{DSmTb2:Dmb:Sec:5:subsec:logtobelfus:eq:4}), it becomes possible to derive $m_i\oplus m_j$ from $m_i$ and $m_j$\,.
Three different cases arise.
\begin{itemize}
\item Axiom m.iii is removed.
Then, the fusion rule~(\ref{dmb:dsmtb2:ch2:rule:eq:2:bis}) of TBM is recovered,
\item Axiom m.iii is used, but $<\Theta>_\Gamma$ verifies the insulation property~(\ref{dmb:dsmtb2:ch2:insulation:eq:1}).
Then, the fusion rule~(\ref{dmb:dsmtb2:ch2:rule:eq:2}) of free DSmT is recovered,
\item Axiom m.iii is used in the general case.
Then, the definition~(\ref{dmb:dsmtb2:ch2:rule:eq:5}) of EMR is recovered.
Moreover, $\Box_i(\phi)\subset\phi$ implies $\mathrm{Bel}_i(\phi)\le p(\phi)$, which is exactly the bound hypothesis~(\ref{dmb:dsmtb2:ch2:Implement:eq:1:1}).
\vspace{7pt}
\end{itemize}
\emph{The logical justification of rule EMR is now completed.}
\section{Conclusion}
\label{DSmTb2:Dmb:Sec:6}
In this chapter, a new fusion rule have been defined for evidence theories.
This rule is computed in order to maximize the entropy of the joint information.
This method provides an adaptive implementation of the independence hypothesis of the sensors.
The rule has been tested on typical examples by means of an algorithmic optimization and by means of a direct computation.
It has been shown that it does not generate conflicts and is compatible with a probabilistic bound interpretation of the belief function.
It is still able to detect truly conflicting sources however, since the optimization may be unfeasible on these cases.
At last, a main contribution of this rule is also that it is derived from an interpretation of evidence theories by means of modal logics.
%
%
%
%********************************
% Bibliographie
%********************************


\begin{thebibliography}{99}
%\bibitem{cox}
%Cox R.~T., \emph{The Algebra of Probale Inference}, Johns Hopkins Press, Baltimore, Maryland, USA, 1961.
%
\bibitem{appriou}
Appriou A.,
\emph{Discrimination multisignal par la th\'eorie de l'\'evidence}
in \emph{D\'ecision et reconnaissance des formes en signal},
Trait\'e IC2, Editions Herm\`es, 2002.
%
\bibitem{blackburn}
Blackburn P., De Rijke M., Venema Y.,
\emph{Modal Logic (Cambridge Tracts in Theoretical Computer Science)},
Cambridge University Press, 2002.
%
\bibitem{DSmTb1:dmb}
Dambreville F.,
\emph{Probabilistic logics related to DSmT}
in \emph{Advances and Applications of DSmT for Information Fusion}, F. Smarandache and J. Dezert editors, American Research Press, Rohoboth, 2004.
%
\bibitem{dempster2}
Dempster A.~P., \emph{A generalization of Bayesian inference}, J. Roy. Statist. Soc., Vol. B, No. 30, pp. 205--247, 1968.
%
\bibitem{dezert}
Dezert J., \emph{Foundations for a new theory of plausible and paradoxical reasoning}, Information \& Security, An international Journal, edited by Prof. Tzv. Semerdjiev, CLPP, Bulg. Acad. of Sciences, Vol. 9, 2002.
%
\bibitem{DezSm:3}
Smarandache F., Dezert J.,
\emph{Information Fusion Based on New Proportional Conflict Redistribution Rules},
8th international conference on information fusion, Philadelphia, 2005.
%
\bibitem{FloJouss}
Florea M.C., Grenier D., Jousselme A.-L., Boss\'e E.,
\emph{Adaptive Combination Rule for Evidence Theory},
accepted at SMCA.
%
\bibitem{martin2}
Osswald C., Martin A.,
\emph{Understanding the large family of Dempster-Shafer theory's fusion operators - a decision-based measure},
9th international conference on information fusion, Florence, 2006.
%
\bibitem{shafer}
Shafer G., \emph{A Mathematical Theory of Evidence}, Princeton Univ. Press, Princeton, NJ, 1976.
%
\bibitem{DSmTBook1}
F. Smarandache and J. Dezert editors, \emph{Advances and Applications of DSmT for Information Fusion}, American Research Press, Rohoboth, 2004.
%
\bibitem{smets}
Smets Ph., \emph{The combination of evidences in the transferable belief model}, IEEE Transactions on Pattern Analysis and Machine Intelligence, Vol. 12, no. 5, pp. 447--458, 1990.
\end{thebibliography}
\end{document}